\documentclass[a4paper,11pt,english]{amsart}

\usepackage[latin1]{inputenc}

\usepackage{calc}
\usepackage{graphicx,palatino,verbatim}

\usepackage{amsmath}

\usepackage[T1]{fontenc}
\usepackage[english]{babel}
\usepackage{newlfont}

\usepackage{amsthm}

\usepackage{amssymb}

\usepackage{amsmath}

\usepackage{eufrak}

\usepackage{latexsym}

\usepackage{young}
\usepackage[mathscr]{eucal}

\usepackage{amsfonts}
\usepackage{syntonly}
\usepackage{graphicx}

\usepackage[letterpaper]{geometry}
\input xy
\xyoption{all}

\newtheorem{teo}{Theorem}

\newtheorem{teorema}{Theorem}[section]

\newtheorem{oss}[teorema]{Remark}

\newtheorem{lemma}[teorema]{Lemma}
\newtheorem{definizione}[teorema]{Definition}

\newtheorem{proposizione}[teorema]{Proposition}

\newtheorem{cor}[teorema]{Corollary}

\begin{document}

\title{Semi-invariants of symmetric quivers of finite type}
\maketitle
\begin{center}
\author{Riccardo Aragona}
\end{center}
\hspace*{2mm}
\begin{center}
Università degli Studi di Roma \textquotedblleft
Tor Vergata\textquotedblright\\
Dipartimento di Matematica\\
Via della Ricerca Scientifica 1, 00133 Rome (Italy)
\end{center}
\hspace*{1mm}
\begin{center}
E-mail: $\quad$ ric$\_$aragona@yahoo.it
\end{center}
\begin{abstract}
\maketitle Let $(Q,\sigma)$ be a symmetric quiver, where
$Q=(Q_0,Q_1)$ is a finite quiver without oriented cycles and
$\sigma$ is a contravariant involution on $Q_0\sqcup Q_1$. The
involution allows us to define a nondegenerate bilinear form $<,>$
on a representation $V$ of $Q$. We shall call the representation
orthogonal if $<,>$ is symmetric and symplectic if $<,>$ is
skew-symmetric. Moreover we can define an action of products of
classical groups on the space of orthogonal representations and on
the space of symplectic representations. For symmetric quivers of
finite type, we prove that the rings of semi-invariants for this
action are spanned by the semi-invariants of determinantal type
$c^V$ and, in the case when matrix defining $c^V$ is
skew-symmetric, by the Pfaffians $pf^V$.
\end{abstract}

\section*{Introduction}
The representations of symmetric quivers, introduced by Derksen
and Weyman in \cite{dw2}, provide a formalization of some problems
related to representations of classical groups.\\
Magyar, Weyman and Zelevinsky in \cite{mwz} classified products of
flag varieties with finitely many orbits under the diagonal action
of general linear groups. The representations of symmetric quivers
could be a tool to solve similar problems for classical groups.\\
Igusa, Orr, Todorov and Weyman in \cite{iotw} generalized the
semi-invariants of quivers to virtual representations of quivers.
They associated, via virtual semi-invariants of quivers, a
simplicial complex $\mathcal{T}(Q)$ with each quiver $Q$. In
particular, if $Q$ is of finite type, then the simplices of
$\mathcal{T}(Q)$ correspond to tilting objects in a corresponding
Cluster category (defined in \cite{bmrrt}). It would be
interesting to carry out a similar construction for
semi-invariants of symmetric quivers of finite type and to relate
it to Cluster algebras (see \cite{fz1} and \cite{fz2}).\\
For usual quivers, Derksen and Weyman in \cite{dw1} proved that
the rings of semi-invariants are spanned by semi-invariants of
determinantal type $c^V$ (introduced by Schofield in \cite{s}),
where $V$ is a representation of the quiver in question. An analogous result has been independently proved with a different approach by Domokos and Zubkov in \cite{dz}.\\ 
In this paper, which is extract from my PhD thesis, supervised by
Professor Jerzy Weyman, whom I gratefully thank, we provide
similar results for symmetric quivers of finite type. In next
paper we will deal with generators of rings of semi-invariants for
symmetric quivers of tame type.\\
A symmetric quiver is a pair $(Q,\sigma)$ where $Q$ is a quiver
(called \textit{underlying quiver of} $(Q,\sigma)$) and $\sigma$
is a contravariant involution on the union of the set of arrows
and the set of vertices of $Q$. The involution allows us to define
a nondegenerate bilinear form $<,>$ on a representation $V$ of
$Q$. We call the pair $(V,<,>)$ orthogonal representation
(respectively symplectic) of $(Q,\sigma)$ if $<,>$ is symmetric
(respectively skew-symmetric). We define $SpRep(Q,\beta)$ and
$ORep(Q,\beta)$ to be respectively the space of symplectic
$\beta$-dimensional representations and the space of orthogonal
$\beta$-dimensional representations  of $(Q,\sigma)$. Moreover we
can define an action of a product of classical groups, which we
call $SSp(Q,\beta)$ in the symplectic case and $SO(Q,\beta)$ in
the orthogonal case, on these space. We describe a set of
generators of the ring of semi-invariants $OSI(Q,\beta)$ of
$ORep(Q,\beta)$ and of the ring of semi-invariants $SpSI(Q,\beta)$
of $SpRep(Q,\alpha)$.\\
Let $(Q,\sigma)$ be a symmetric quiver and $V$ a representation of
the underlying quiver $Q$ such that
  $\langle\underline{dim}\,V,\beta\rangle=0$, where $\langle\cdot,\cdot\rangle$ is the Euler form of $Q$. Let
$$
0\longrightarrow P_1\stackrel{d^V}{\longrightarrow}
P_0\longrightarrow V\longrightarrow 0
$$
be the canonical projective resolution of $V$ (see \cite{r1}). We
define the semi-invariant $c^V:=det(Hom_Q(d^V,\cdot))$ of $
OSI(Q,\beta)$ and $ SpSI(Q,\beta)$ (see
  \cite{dw1} and \cite{s}).\\
   Let $C^+$ be the Coxeter functor and let $\nabla$ be the duality functor.
    We will
  prove in the symmetric case the following
  \begin{teo}\label{p1}
Let $(Q,\sigma)$ be a symmetric quiver of finite type and let
$\beta$ be a symmetric dimension vector. The ring $SpSI(Q,\beta)$
is generated by semi-invariants
\begin{itemize}
\item[(i)] $c^V$ if $V\in Rep(Q)$ is such that
$\langle\underline{dim}\,V,\beta\rangle=0$,
\item[(ii)] $pf^V$ if $V\in Rep(Q)$ is such that
$\langle\underline{dim}\,V,\beta\rangle=0$, $C^+ V=\nabla V$ and
the almost split sequence $0\rightarrow\nabla V\rightarrow
Z\rightarrow V\rightarrow 0$ has the middle term $Z$ in $ORep(Q)$.
\end{itemize}
\end{teo}
\begin{teo}\label{p2}
Let $(Q,\sigma)$ be a symmetric quiver of finite type and let
$\beta$ be a symmetric dimension vector. The ring $OSI(Q,\beta)$
is generated by semi-invariants
\begin{itemize}
\item[(i)] $c^V$ if $V\in Rep(Q)$ is such that
$\langle\underline{dim}\,V,\beta\rangle=0$,
\item[(ii)] $pf^V$ if $V\in Rep(Q)$ is such that
$\langle\underline{dim}\,V,\beta\rangle=0$, $C^+ V=\nabla V$ and
the almost split sequence $0\rightarrow\nabla V\rightarrow
Z\rightarrow V\rightarrow 0$ has the middle term $Z$ in
$SpRep(Q)$.
\end{itemize}
\end{teo}
A similar result has been obtained by Lopatin in \cite{l} in a different setting, using ideas from \cite{lz}.\\
The strategy of the proofs is the following. First we set the
technique of reflection functors on the symmetric quivers. Then we
prove that we can reduce Theorems \ref{p1} and \ref{p2}, by this
technique, to particular orientations of the symmetric quivers.
Finally, we check Theorems \ref{p1} and \ref{p2} for these
orientations.\\
In the first section we give general notions and results about
usual and symmetric quivers and their representations. We state
main results \ref{p1} and \ref{p2} and we recall some results of
representations of general linear groups and of invariant
theory.\\
In the second section we adjust to symmetric quivers the technique
of reflection functors and we prove general results about
semi-invariants of symmetric quivers.\\
In the third section we check we can reduce Theorems \ref{p1} and
\ref{p2} to a particular orientation of symmetric quivers of
finite type. Finally, using classical invariant theory and the
technique of Schur functors, we prove case by case Theorems
\ref{p1} and \ref{p2} for symmetric quivers of finite type with
this orientation.

\section*{Acknowledgments}
I am deeply grateful to my PhD thesis advisor Jerzy Weyman for his constant guidance and for all suggestions and ideas he shared me. I am
grateful to Elisabetta Strickland for the great
helpfulness she showed me during the three years of my graduate studies. I also
wish to thank Fabio Gavarini and Alessandro
D'Andrea for several useful discussions. I thank Alexander Zubkov for pointing out some interesting references.

\section{Preliminary results}
Throughout the paper $\Bbbk$ denotes an algebraically closed field
of characteristic 0.
 \subsection{Representations of quivers}
A quiver $Q$ is a quadruple $(Q_0,Q_1,t,h)$ where $Q_0$ is a set
of vertices, $Q_1$ is a set of arrows and $t,h:Q_1\rightarrow Q_0$
are two maps which assign to each arrow $a\in Q_1$ respectively
its tail $ta\in Q_0$ and its head $ha\in Q_0$. Throughout the
paper we consider quivers $Q$ without oriented cycles, i.e. in
which there are no paths
$a_1\cdots a_n$ such that $ta_1=ha_n$.\\
A representation $V$ of $Q$ is a family of finite dimensional
vector spaces
 $\{V(x)|x\in Q _0\}$ and of linear maps
$\{V(a):V(ta)\rightarrow V(ha)\}_{a\in Q_1}$. The dimension vector
of $V$ is a function $\underline{dim}(V):Q_0\rightarrow\mathbb{N}$
defined by $\underline{dim}(V)(x):=dim V(x)$. For a dimension
vector $\alpha$ we have
$$
Rep(Q,\alpha)=\bigoplus_{a\in
Q_0}Hom(\Bbbk^{\alpha(ta)},\Bbbk^{\alpha(ha)})
$$
the variety of representations of $Q$ of dimension $\alpha$.
Moreover we define
$$
SL(Q,\alpha)=\prod_{x\in Q_0}SL(\alpha(x))<\prod_{x\in
Q_0}GL(\alpha(x))=GL(Q,\alpha).
$$
The action of these groups on $Rep(Q,\alpha)$ is defined by
$$
g\cdot V=\{g_{ha}V(a)g_{ta^{-1}}\}_{a\in Q_1}
$$
where $g=(g_x)_{x\in Q_0}\in GL(Q,\alpha)$ and $V\in
Rep(Q,\alpha)$.\\
A morphism $f:V\rightarrow W$ of two
representations is a family of linear maps $\{f(x):V(x)\rightarrow
W(x)|\,f(ha)V(a)=W(a)f(ta)\forall a\in Q_1\}_{x\in Q_0}$. We
denote the space of morphisms from $V$ to $W$ by $Hom_Q(V,W)$ and
the space of extensions of $V$ by $W$ by $Ext_Q^1(V,W)$.
\begin{definizione}\label{fe}
The non symmetric bilinear form on the space $\mathbb{Z}^{Q_0}$ of
dimension vectors given by
$$
\langle\alpha,\beta\rangle=\sum_{x\in
Q_0}\alpha(x)\beta(x)-\sum_{a\in Q_1}\alpha(ta)\beta(ha)
$$
is the Euler form of $Q$, where $\alpha,\beta\in\mathbb{Z}^{Q_0}$.
\end{definizione}
A vertex $x\in Q_0$ is said to be a sink (resp. a source) of $Q$
if
$x=ha$ (resp. $x=ta$) for every $a\in Q_1$ connected to $x$.\\
Let $x\in Q_0$ be a sink (resp. a source) of a quiver $Q$ and let
$\{a_1,\ldots,a_k\}$ be the arrows connected to $x$. We define the
quiver $c_x(Q)$ as follows
$$c_x(Q)_0=Q_0$$
$$c_x(Q)_1=\{c_x(a)|a\in Q_1\}$$
where $tc_x(a_i)=ha_i$,   $hc_x(a_i)=ta_i$ for every
$i\in\{1,\ldots,k\}$ and $tc_x(b)=tb$, $hc_x(b)=hb$ for every
$b\in Q_1\setminus\{a_1,\ldots,a_k\}$. Moreover, we can define the
reflection $c_x:\mathbb{Z}^{Q_0}\rightarrow\mathbb{Z}^{Q_0}$ given
for $\alpha\in\mathbb{Z}^{Q_0}$ by formula
$$
c_x(\alpha)(y)=\left\{\begin{array}{ll}
\alpha(y) & if\quad y\ne x\\
\sum_{i=1}^k \alpha(ta_i)-\alpha(x) & \textrm{otherwise}.
\end{array}\right.
$$
Finally it is known (see \cite{bgp} and \cite{dr}) that for every
$x\in Q_0$ sink or a source of $Q$ we can define respectively the
functors
$$
C_x^+:Rep(Q,\alpha)\rightarrow Rep(c_x(Q),c_x(\alpha))
$$
and
$$C_x^-:Rep(Q,\alpha)\rightarrow Rep(c_x(Q),c_x(\alpha))$$
called \textit{reflection functors}.
\begin{definizione}
Let $Q$ be a quiver with $n$ vertices without oriented cycles. We
choose the numbering $(x_1,\ldots,x_n)$ of vertices such that
$ta>ha$ for every $a\in Q_1$. We define
$$
C^+:=C^+_{x_n}\cdots C^+_{x_1}\quad\textrm{and}\quad
C^-:=C^-_{x_1}\cdots C^-_{x_n}.
$$
The functors $C^+,C^-:Rep(Q)\rightarrow Rep(Q)$ are called Coxeter
functors. \end{definizione} One proves that these functors don't
depend on the choice of numbering of vertices (see \cite{ass}
chap. VII Lemma 5.8).

 \subsection{Symmetric quivers}\begin{definizione}
A symmetric quiver is a pair $(Q,\sigma)$ where $Q$ is a quiver
(called underlying quiver of $(Q,\sigma)$) and $\sigma$ is an
involution on $Q_0\sqcup Q_1$ such that
\begin{itemize}
   \item[(i)] $\sigma(Q_0)=Q_0$ and $\sigma(Q_1)=Q_1$,
   \item[(ii)] $t\sigma(a)=\sigma(ha)$ and $h\sigma(a)=\sigma(ta)$ for all $a\in
   Q_1$,
   \item[(iii)] $\sigma(a)=a$ whenever $a\in Q_1$ and $\sigma(ta)=ha$.
   \end{itemize}
 \end{definizione}
 Let $V$ be a representation of the underlying quiver $Q$ of a
 symmetric quiver $(Q,\sigma)$. We define the duality functor
 $\nabla:Rep(Q)\rightarrow Rep(Q)$ such that $\nabla
 V(x)=V(\sigma(x))^*$ for every $x\in Q_0$ and $\nabla
 V(a)=-V(\sigma(a))^*$ for every $a\in Q_1$. If $f:V\rightarrow W$
 is a morphism of representations $V,W\in Rep(Q)$, then $\nabla
 f:\nabla W\rightarrow \nabla V$ is defined by $\nabla
 f(x)=f(\sigma(x))^*$, for every $x\in Q_0$. We call $V$
 \textit{selfdual}
 if $\nabla V=V$.
  \begin{definizione}\label{defortsp}
 An orthogonal (resp. symplectic) representation of a symmetric quiver $(Q,\sigma)$ is a pair $(V,<\cdot,\cdot >)$, where $V$ is
  a representation of the underlying quiver $Q$ with a nondegenerate symmetric (resp. skew-symmetric) scalar product $<\cdot,\cdot >$ on $\bigoplus_{x\in Q_0}V(x)$
   such that
   \begin{itemize}
   \item[(i)] the restriction of $<\cdot,\cdot>$ to $V(x)\times V(y)$ is 0 if $y\neq\sigma(x)$,
   \item[(ii)] $<V(a)(v),w>+<v,V(\sigma(a))(w)>=0$ for all $v\in V(ta)$ and all $w\in V(\sigma(a))$.
   \end{itemize}
  \end{definizione}
  By properties \textit{(i)} and \textit{(ii)} of definition \ref{defortsp}, an orthogonal or symplectic representation\\ $(V,<\cdot,\cdot>)$ of a symmetric quiver is
  selfdual.\\
 We shall say that a dimension vector $\alpha$ is \textit{symmetric} if
 $\alpha(x)=\alpha(\sigma(x))$ for every $x\in Q_0$. Since each
 orthogonal or symplectic representation is selfdual, then
 dimension vector of an orthogonal (resp. symplectic)
 representation, which we shall call respectively orthogonal and symplectic dimension vector, is symmetric.
\begin{definizione}\label{defindspo}
An orthogonal (respectively symplectic) representation is called
indecomposable orthogonal (respectively indecomposable symplectic)
if it cannot be expressed as a direct sum of orthogonal
(respectively symplectic) representations.
\end{definizione}
 \begin{definizione}
A symmetric quiver is said to be of finite representation type if
it has only finitely many indecomposable orthogonal (resp.
symplectic) representations up to isomorphisms.
\end{definizione}
Derksen and Weyman classified the symmetric quiver of finite type
in \cite{dw2}
\begin{teorema}\label{ctf}
A symmetric quiver $(Q,\sigma)$ is  of finite type if and only if
the underlying quiver $Q$ is of type $A_n$.
\end{teorema}
\begin{proof} See \cite{dw2}, theorem 3.1 and proposition 3.3.\end{proof}
We describe the space of orthogonal (resp. symplectic)
representations of a symmetric quiver $(Q,\sigma)$.\\
We denote
 $ Q^{\sigma}_0$
 (respectively $Q^{\sigma}_1$) the set of vertices (respectively
arrows) fixed by $\sigma$. Thus we have partitions
$$
Q_0 = Q^+_ 0 \sqcup Q^{\sigma}_ 0 \sqcup Q^-_0
$$
$$
 Q_1 = Q^+_ 1 \sqcup Q^{\sigma}_ 1 \sqcup
Q^-_1
$$
 such that $Q^-_0 = \sigma(Q^+_ 0 )$ and $Q^-_1 = \sigma(Q^+_ 1 )$, satisfying:
 \begin{itemize}
\item[i)] $\forall a \in Q^+_ 1$ , either $\{ta, ha\} \subset Q^+_ 0$ or one of the elements in $\{ta, ha\}$ is
in $Q^+_ 0$ while the other is in $Q^{\sigma}_ 0$;
\item[ii)] $\forall x\in Q^+_0$, if $a\in Q_1$ with $ta = x$ or $ha = x$, then $a\in Q^+_ 1 \sqcup
Q^{\sigma}_1$.
\end{itemize}
The space of orthogonal $\alpha$-dimensional
  representations of a symmetric quiver $(Q,\sigma)$ can be
  identified with
  \begin{equation}
   ORep(Q,\alpha)=\bigoplus_{a\in Q^{+}_1 }Hom(\Bbbk^{\alpha(ta)},\Bbbk^{\alpha(ha)})\oplus\bigoplus_{a\in
Q^{\sigma}_1}\bigwedge^2(\Bbbk^{\alpha(ta)})^*.
  \end{equation}
  The space of symplectic $\alpha$-dimensional representations
  can be identified with
 \begin{equation}
 SpRep(Q,\alpha)=\bigoplus_{a\in Q^+_1}Hom(\Bbbk^{\alpha(ta)},\Bbbk^{\alpha(ha)})\oplus\bigoplus_{a\in Q^{\sigma}_1 }S_2(\Bbbk^{\alpha(ta)})^*.
   \end{equation}
  We define the group
\begin{equation}
 SO(Q,\alpha)=\prod_{x\in
Q^+_0}SL(\alpha(x))\times\prod_{x\in Q^{\sigma}_0}SO(\alpha(x)),
   \end{equation}
  where $SO(\alpha(x))$ is the group of special orthogonal
  transformations for the symmetric form
  $<\cdot,\cdot>$ restricted to $V(x)$.\\
  \\
  Assuming that $\alpha(x)$ is even
  for every $x\in Q_0^{\sigma}$, we define the group
\begin{equation}
   SSp(Q,\alpha)=\prod_{x\in Q^+_0}SL(\alpha(x))\times\prod_{x\in
   Q^{\sigma}_0}Sp(\alpha(x)),
   \end{equation}
where $Sp(\alpha(x))$ is the group of isometric
  transformations for the skew-symmetric form\\
  $<\cdot,\cdot>$ restricted to $V(x)$.\\
  \\
The action of these groups is defined by
$$
g\cdot V=\{g_{ha}V(a)g_{ta}^{-1}\}_{a\in Q_1^+\cup Q_1^{\delta}}
$$
where $g=(g_x)_{x\in Q_0}\in SO(Q,\alpha)$ (respectively
  $g\in SSp(Q,\alpha)$) and $V\in ORep(Q,\alpha)$ (respectively in $SpRep(Q,\alpha)$). In particular we can suppose $g_{\sigma(x)}=(g_x^{-1})^t$
  for every $x\in Q_0$.
  \subsection{Semi-invariants of quivers without oriented cycles and main
  results}In this section first we define semi-invariants which
  appear in main results of paper and we describe some property of
  these for any quiver $Q$; then we state main theorems.\\
  Let $Q$ be a quiver with $n$ vertices. We denote
$$
SI(Q,\alpha)=\Bbbk[Rep(Q,\alpha)]^{SL(Q,\alpha)}
$$
the ring of semi-invariants of a quiver $Q$.\\
For every $g\in GL(Q,\alpha)$ the character $\tau$ at $g$ is
$\tau(g)=det(g)^{\chi_1}\cdots det(g)^{\chi_n}$, where
$\chi=(\chi_1,\ldots,\chi_n)\in \mathbb{Z}^n$ is also called
\textit{weight} if $\tau$ is a weight for some semi-invariant. So
the ring $SI(Q,\alpha)$ decomposes in graded components as
  $$
  SI(Q,\alpha)=\bigoplus_{\tau\in char(GL(Q,\alpha))}SI(Q,\alpha)_{\tau}
  $$
 where $SI(Q,\alpha)_{\tau} =\big\{f\in\Bbbk[Rep(Q,\alpha)]|g\cdot f=\tau(g)f\;\forall g\in
 GL(Q,\alpha)\big\}$.\\
We define the semi-invariants which appear in the main theorems.\\
For every $V\in Rep(Q,\alpha)$, we can construct a projective
resolution, called \textit{canonical resolution of} $V$:
\begin{equation}\label{Rr}
0\longrightarrow\bigoplus_{a\in Q_1}V(ta)\otimes
P_{ha}\stackrel{d^V}{\longrightarrow}\bigoplus_{x\in
Q_0}V(x)\otimes
P_{x}\stackrel{p_V}{\longrightarrow}V\longrightarrow 0
\end{equation}
where $P_x$ is the indecomposable projective associated to vertex
$x$ for every $x\in Q_0$ (see \cite{ass}), $d^V|_{V(ta)\otimes
P_{ha}}(v\otimes e_{ha})=V(a)(v)\otimes e_{ha}-v\otimes a$ and
$p_V|_{V(x)\otimes P_{x}}(v)=v\otimes e_{x}$. Applying the functor
$Hom_Q(\cdot,W)$ to $d^V$ for $W\in Rep(Q,\beta)$, we have that
the matrix associated to $Hom_Q(d^V,W)$ is square if and only if
$\langle\alpha,\beta\rangle=0$ (see \cite{s} Lemma 1.2).
\begin{definizione}
For $V\in Rep(Q,\alpha)$ such that $\langle\alpha,\beta\rangle=0$,
where $\beta\in\mathbb{N}^n$, we define
$$
\begin{array}{rcl}
c^V:Rep(Q,\beta)&\longrightarrow&\Bbbk\\
W&\longmapsto&c^V(W)=det(Hom_Q(d^V,W)).
\end{array}
$$
These are semi-invariants of weight $\langle\alpha,\cdot\rangle$,
called Schofield semi-invariants (see \cite{s} Lemma 1.4).
\end{definizione}
\begin{oss}\label{qp}
\begin{itemize}
\item[(i)] Any projective resolution of $V$ can be used to
calculate $c^V$ (see \cite{s}). Moreover if $P$ is a projective
representation, then $c^P=0$.
\item[(ii)] If
$\langle\underline{dim}V,\underline{dim}W\rangle=0$, then
$c^V(W)=0$ if and only if $Hom_Q(V,W)\neq 0$ (see \cite{dw1}).
\end{itemize}
\end{oss}
Now we formulate the result of Derksen and Weyman about the set of
generators of the ring of semi-invariants of a quiver without
oriented cycles $Q$.
\begin{teorema}[Derksen-Weyman]\label{dw}
Let $Q$ be a quiver without oriented cycles and let $\beta$ be a
dimension vector. The ring $SI(Q,\beta)$ is spanned by
semi-invariants of the form $c^V$ of weight $\langle
\underline{dim}(V),\cdot\rangle$, for which $\langle
\underline{dim}(V),\beta\rangle=0$.
\end{teorema}
\begin{proof}
See \cite{dw1} Theorem 1.
\end{proof}
We give some property of Schofield semi-invariants.
\begin{lemma}\label{cVcV'}
Suppose that $V'$, $V$, $V''$ and $W$ are representations of $Q$,
that $\langle\underline{dim}(V),\underline{dim}(W)\rangle=0$ and
that there are exact sequences
$$
0\rightarrow V'\rightarrow V\rightarrow V''\rightarrow 0
$$
then
\begin{itemize}
\item[(i)] If $\langle\underline{dim}(V'),\underline{dim}(W)\rangle<0$, then $c^V(W)=0$
\item[(ii)] If $\langle\underline{dim}(V'),\underline{dim}(W)\rangle=0$, then $c^{V}(W)=c^{V'}(W)
c^{V''}(W)$.
\end{itemize}
\end{lemma}
\begin{proof} See \cite{dw1} Lemma 1.\end{proof}
We recall definition and properties of the
  \textit{Pfaffian} of a skew-symmetric matrix.\\
  Let $A=(a_{ij})_{1\leq i,j\leq 2n}$ be a skew-symmetric $2n\times 2n$
  matrix. Given $2n$ vectors $x_1,\ldots,x_{2n}$ in $\Bbbk^{2n}$,  we define
  $$
  F_A(x_1,\ldots,x_{2n})=\sum_{{i_1<j_1,\ldots,i_n<j_n \atop i_1<\ldots<i_n}}sgn(s)\prod_{i=1}^n(x_{s(2i-1)},x_{s(2i)}),
  $$
  where $sgn(s)$
  is the sign of permutation
  $$
  s=\left[\begin{array}{ccccc} 1 & 2
  &\ldots & 2n-1 & 2n \\
  i_1 & j_1 & \ldots & i_n & j_n \end{array}\right]
  $$
  and $(\cdot,\cdot)$ is the
  skew-symmetric bilinear form associated to $A$. So $F_A$ is a
  skew-symmetric multilinear function of $x_1,\ldots,x_{2n}$.
  Since, up to a scalar, the only one skew-symmetric
  multilinear function of $2n$ vectors in $\Bbbk^{2n}$ is the
  determinant, there is a complex number $Pf(A)$, called
  \textit{Pfaffian of} $A$, such that
  $$
  F_A(x_1,\ldots,x_{2n})=Pf(A)det[x_1,\ldots,x_{2n}]
  $$
  where $[x_1,\ldots,x_{2n}]$ is the matrix which has the vector
  $x_i$ for $i$-th column. In particular we note that
  $$
 Pf(A)=\sum_{{i_1<j_1,\ldots,i_n<j_n \atop i_1<\ldots<i_n}}sgn\left(\left[\begin{array}{ccccc} 1 & 2
  &\ldots & 2n-1 & 2n \\
  i_1 & j_1 & \ldots & i_n & j_n \end{array}\right]\right)a_{1_1j_1}\cdots
  a_{i_nj_n}.
  $$
  \begin{proposizione}
  Let $A$ be a skew-symmetric $2n\times 2n$ matrix.
  \begin{itemize}
  \item[(i)] For every invertible $2n\times 2n$ matrix $B$,
  $$
  Pf(BAB^t)=det(B)Pf(A);
  $$
  \item[(ii)] $det(A)=Pf(A)^2$.
  \end{itemize}
  \end{proposizione}
  \begin{proof} See \cite{p}, chap. 5 sec. 3.6. \end{proof}
Let $V\in Rep(Q,\alpha)$ and let $\beta$ be a dimension vector
such that $\langle\alpha,\beta\rangle=0$ and
$Hom_Q(d^V_{min},\cdot)$ is skew-symmetric on $Rep(Q,\beta)$, we
can define
$$
\begin{array}{rcl}
pf^V:Rep(Q,\beta)&\longrightarrow&\Bbbk\\
W&\longmapsto&pf^V(W)=Pf(Hom_Q(d^V_{min},W)).
\end{array}
$$
In this work we describe a set of generators of the rings of
semi-invariants of symmetric quivers of finite type.\\
Let $\alpha$ be a dimension vector of an orthogonal or symplectic
representation, we denote
$$
OSI(Q,\alpha):=\Bbbk[ORep(Q,\alpha)]^{SO(Q,\alpha)}\quad\textrm{and}\quad
SpSI(Q,\alpha):=\Bbbk[SpRep(Q,\alpha)]^{SSp(Q,\alpha)}
$$
respectively the ring of orthogonal semi-invariants and the ring
of symplectic semi-invariants of a symmetric quiver
$(Q,\sigma)$.\\
We state the main theorems
\begin{teorema}\label{tp1}
Let $(Q,\sigma)$ be a symmetric quiver of finite type and let
$\beta$ be a symmetric dimension vector. The ring $SpSI(Q,\beta)$
is generated by semi-invariants
\begin{itemize}
\item[(i)] $c^V$ if $V\in Rep(Q)$ is such that
$\langle\underline{dim}\,V,\beta\rangle=0$,
\item[(ii)] $pf^V$ if $V\in Rep(Q)$ is such that
$\langle\underline{dim}\,V,\beta\rangle=0$, $C^+ V=\nabla V$ and
the almost split sequence $0\rightarrow\nabla V\rightarrow
Z\rightarrow V\rightarrow 0$ has the middle term $Z$ in $ORep(Q)$.
\end{itemize}
\end{teorema}
\begin{teorema}\label{tp2}
Let $(Q,\sigma)$ be a symmetric quiver of finite type let $\beta$
be a symmetric dimension vector. The ring $OSI(Q,\beta)$ is
generated by semi-invariants
\begin{itemize}
\item[(i)] $c^V$ if $V\in Rep(Q)$ is such that
$\langle\underline{dim}\,V,\beta\rangle=0$,
\item[(ii)] $pf^V$ if $V\in Rep(Q)$ is such that
$\langle\underline{dim}\,V,\beta\rangle=0$, $C^+ V=\nabla V$ and
the almost split sequence $0\rightarrow\nabla V\rightarrow
Z\rightarrow V\rightarrow 0$ has the middle term $Z$ in
$SpRep(Q)$.
\end{itemize}
\end{teorema}
The strategy of the proofs is the following. First we set the
technique of reflection functors on the symmetric quivers. Then we
prove that we can reduce theorems \ref{tp1} and \ref{tp2}, by this
technique, to particular orientations of the symmetric quivers.
Finally, we check theorems \ref{tp1} and \ref{tp2} for these
orientations using technique of Schur functors.
\begin{definizione}\label{popssp}
Let $(Q,\sigma)$ be a symmetric quiver. We will say that $V\in
Rep(Q)$ satisfies property \textit{(Op)} if
\begin{itemize}
\item[(i)] $V=C^-\nabla V$
\item[(ii)] the almost split sequence $0\rightarrow\nabla V\rightarrow
Z\rightarrow V\rightarrow 0$ has the middle term $Z$ in $ORep(Q)$.
\end{itemize}
Similarly we will say that $V\in Rep(Q)$ satisfies property
\textit{(Spp)} if
\begin{itemize}
\item[(i)] $V=C^-\nabla V$
\item[(ii)] the almost split sequence $0\rightarrow\nabla V\rightarrow
Z\rightarrow V\rightarrow 0$ has the middle term $Z$ in
$SpRep(Q)$.
\end{itemize}
\end{definizione}

\subsection{Invariant theory and Schur modules}
Let $G$ be an algebraic group, $V$ a rational representation of
$G$ and $\Bbbk[V]$ the algebra of regular functions of $V$.\\
If $\mathcal{X}(G)$ is the set of characters of $G$, then the ring
of the semi-invariants of $G$ on $V$ is defined by
$$
SI(G,V)=\bigoplus_{\chi\in\mathcal{X}(G)}SI(G,V)_\chi
$$
where $SI(G,V)_\chi=\{f\in\Bbbk[V]|g\cdot f=\chi(g)f,\,\forall
g\in G\}$ is called weight space of weight $\chi$. The following
lemma describes $SI(G,V)$ in the case when $G$ has an open orbit
on $V$.
\begin{lemma}[Sato-Kimura]\label{sk}
Let $G$ be a connected linear algebraic group and $V$ a rational
representation of $G$. We suppose that the action of $G$ on $V$
has an open orbit. Then $SI(G,V)$ is a polynomial $\Bbbk$-algebra
and the weights of the generators of $SI(G,V)$ are linearly
independent in $\mathcal{X}(G)$. Moreover, the dimensions of the
spaces $SI(G,V)_\chi$ are 0 or 1.
\end{lemma}
\begin{proof}
See \cite{sk} sect. 4, Lemma 4 and Proposition 5.
\end{proof}
Let $G=GL_n(\Bbbk)$ be the general linear group over $\Bbbk$.
There exists an isomorphism $\mathbb{Z}\cong \mathcal{X}(G)$ which
sends an element $a$ of $\mathbb{Z}$ in $(det)^a$ (where $det$
associates to $g\in G$ its determinant). We identify $G$ with the
group $GL(V)$ of linear automorphisms of a vector space $V$ of
dimension $n$. So we have
 $$
 SI(G,V)=\Bbbk[V]^{SL(V)}.
 $$
Let $T$ and $\mathcal{X}(T)$ respectively be the maximal torus in
$G$ (i.e. the group of diagonal matrices) and the set of
characters of $T$. The irreducible rational representations of $G$
are parametrized by the set
$$
\mathcal{X}^+(T)=\{\lambda=(\lambda_1,\ldots,\lambda_n)\in\mathbb{Z}^n|\lambda_1\geq\cdots\geq\lambda_n\}
$$
of the integral dominant weights for $GL_n(\Bbbk)$. The
irreducible rational representations $S_{\lambda}V$ of
$G=GL_n(\Bbbk)$ corresponding to the dominant weight
$\lambda\in\mathcal{X}^+(T)$ are called Schur modules. In the case
when $\lambda_n\geq 0$ (i.e. $\lambda$ is a partition of
$\lambda_1+\cdots+\lambda_n$), a description of $S_{\lambda}V$ is
given in \cite{abw}, \cite{f} and \cite{p}. For every
$\lambda\in\mathcal{X}^+(T)$, we can define $S_{\lambda}V$ as
follows
$$
S_{(\lambda_1,\ldots,\lambda_n)}V=S_{(\lambda_1-\lambda_n,\ldots,\lambda_{n-1}-\lambda_n,0)}V\otimes(\bigwedge^n
V)^{\otimes\lambda_n}.
$$
Let $\lambda=(\lambda_1,\ldots,\lambda_n)$ be a partition of
$|\lambda|:=\lambda_1+\cdots+\lambda_n$. We call \textit{height}
of $\lambda$, denoted by $ht(\lambda)$, the number $k$ of nonzero
components of $\lambda$ and we denote the transpose of $\lambda$
by $\lambda'$.
\begin{teorema}[Properties of Schur modules]\label{pms}
Let $V$ be vector space of dimension $n$ and $\lambda$ be an
integral dominant weight.
\begin{itemize}
\item[(i)] $S_{\lambda}V=0\Leftrightarrow ht(\lambda)>0$.
\item[(ii)] $dim\,S_{\lambda}V=1\Leftrightarrow
\lambda=(\overbrace{k,\ldots,k}^n)$ for some $n\in\mathbb{Z}$.
\item[(iii)] $\left(S_{(\lambda_1,\ldots,\lambda_n)}V\right)^*\cong
S_{(\lambda_1,\ldots,\lambda_n)}V^*\cong
S_{(-\lambda_n,\ldots,-\lambda_1)}V.$
\end{itemize}
\end{teorema}
\begin{proof}
See Theorem 6.3 in \cite{fh}.
\end{proof}
\begin{teorema}[Cauchy formulas]\label{fc}
Let $V$ and $W$ be two finite dimensional vector spaces. Then
\begin{itemize}
\item[(i)] As representations of $GL(V)\times GL(W)$,
$$
S_d(V\otimes W)=\bigoplus_{|\lambda|=d}S_{\lambda}V\otimes
S_{\lambda}W\quad\textrm{and}\quad \bigwedge^d(V\otimes
W)=\bigoplus_{|\lambda|=d}S_{\lambda}V\otimes S_{\lambda'}W.
$$
\item[(ii)] As representations of
$GL(V)$,
$$S_d(S_2(V))=\bigoplus_{|\lambda|=d}S_{2\lambda}V\quad\textrm{and}\quad S_d(\bigwedge^2(V))
=\bigoplus_{|\lambda|=d}S_{2\lambda'}V,
$$
where $2\lambda=(2\lambda_1,\ldots,2\lambda_k)$ if
$\lambda=(\lambda_1,\ldots,\lambda_k)$.
\end{itemize}
\end{teorema}
\begin{proof}
See \cite{p} chap. 9 sec. 6.3 and sec 8.4, chap 11 sec. 4.5.
\end{proof}
The decomposition of tensor product of Schur modules è
$$
S_{\lambda}V\otimes S_{\mu}V=\bigoplus_{\nu}c^{\nu}_{\lambda\mu}
S_{\nu}V,
$$
where the coefficients $c^{\nu}_{\lambda\mu}$ are called
Littlewood-Richardson coefficients.\\
There exists a combinatorial
formula to calculate $c^{\nu}_{\lambda\mu}$, called
\textit{Littlewood-Richardson rule}
(see in \cite{p} chap. 12 sec. 5.3).\\
 Finally we state other two results on Schur modules and invariant theory.
 \begin{proposizione}\label{i1}
 Let $V$ be a vector space of dimension $n$.
 $$
 (S_{\lambda}V)^{SL(V)}\neq 0\Longleftrightarrow\lambda=(k^n)
 $$
 for some $k$ and in this case $S_{\lambda}V$, and so also $(S_{\lambda}V)^{SL(V)}$, have dimension one.
 \end{proposizione}
 \begin{proof}
See Corollary p. 388 in \cite{p}.
\end{proof}
 \begin{proposizione}\label{i2}
  Let $V$ be a vector space of dimension $n$ and let $\lambda$ and $\mu$ be two integral dominant weights. Then\\
  $$
  (S_{\lambda}V\otimes S_{\mu}V)^{SL\,V}\ne 0
  $$
  $$
  \Longleftrightarrow
 $$
$$
 \lambda_i-\lambda_{i+1}=\mu_{n-i}-\mu_{n-i+1}
  $$
  for every $i\in\{1,\ldots,n-1\}$ and in this case the semi-invariant is unique (up to a non zero scalar) and has weight
  $\lambda_1+\mu_n=\lambda_2+\mu_{n-1}=\cdots=\lambda_n+\mu_1$.
  \end{proposizione}

  \begin{proof}It is a corollary of (7.11) in \cite{m} chap. 1 sec. 5.\end{proof}
Let $Sp(2n)=\{A\in GL_{2n}(\Bbbk)|AJA=J\}$ be
  the simplectic group, let $O(n)=\{A\in GL_n(\Bbbk)|A^tA=I\}$
  be the orthogonal group and $SO(n)=\{A\in O(n)|det\,A=1\}$ be the special orthogonal group,
   where I is the identity matrix and
  $J=\left(\begin{array}{cc} 0 & I\\
  -I & 0 \end{array}\right)$.
  \begin{proposizione}\label{i3}
  Let $V$ be an orthogonal space of dimension $n$ and let $W$ be a
  symplectic space of dimension $2n$.
  \begin{itemize}
  \item[(a)] $ dim\,(S_{\lambda}V)^{O(V)}=\left\{\begin{array}{ll}1 &
  \textrm{if}\quad \lambda=2\mu\\
  0 & \textrm{otherwise}
  \end{array}\right.$,
\item[(b)] $ dim\,(S_{\lambda}V)^{SO(V)}=\left\{\begin{array}{ll}1 &
  \textrm{if}\quad \lambda=2\mu+(k^n)\\
  0 & \textrm{otherwise}\end{array}\right.$,
  \item[(c)] $ dim\,(S_{\lambda}W)^{Sp(W)}=\left\{\begin{array}{ll}1 &
  \textrm{if}\quad \lambda=2\mu'\\
  0 & \textrm{otherwise}\end{array}\right.$
  \end{itemize}
  for some partition $\mu$ and for some $k\in\mathbb{N}$.
  \end{proposizione}
  \begin{proof} See \cite{p} chap. 11 cor. 5.2.1 and 5.2.2.
  \end{proof}

\section{Reflection functors and semi-invariants of symmetric quivers}
\subsection{Reflection functors for symmetric quivers}
We adjust the technique of reflection functors to symmetric
quivers.
\begin{definizione}
Let $(Q,\sigma)$ be a symmetric quiver. A sink (resp. a source)
$x\in Q_0$ is called admissible if there are no arrows connecting
$x$ and $\sigma(x)$.
\end{definizione}
By definition of $\sigma$, $x$ is a sink (resp. a source) if and
only if $\sigma(x)$ is a source (resp. a sink). We call
$(x,\sigma(x))$ \textit{the admissible sink-source pair}. So we
can define $c_{(x,\sigma(x))}:=c_{\sigma(x)}c_x$.
\begin{lemma}
If $(Q,\sigma)$ is a symmetric quiver and $x$ is an admissible
sink or source, then $(c_{(x,\sigma(x))}(Q),\sigma)$ is symmetric.
\end{lemma}
\begin{proof}
It follows from definition of $\sigma$.
\end{proof}
\begin{definizione}
Let $(Q,\sigma)$ be a symmetric quiver. A sequence
$x_1,\ldots,x_m\in Q_0$ is an admissible sequence of sinks (or
sources) for admissible sink-source pairs if $x_{i+1}$ is an
admissible sink (resp. source) in $c_{(x_i,\sigma(x_i))}\cdots
c_{(x_1,\sigma(x_1))}(Q)$ for $i=1,\ldots,m-1$.
\end{definizione}
One proves by a simple combinatorial argument the following
\begin{proposizione}\label{Q=Q'}
Let $(Q,\sigma)$ and $(Q',\sigma)$ be two symmetric connected
quivers without cycles, with the same underlying graph but with
different orientations. Then there exists an admissible sequence
of sinks (or source) for admissible sink-source pairs
$x_1,\ldots,x_m\in Q_0$ such that
$$
Q'=c_{(x_m,\sigma(x_m))}\cdots c_{(x_1,\sigma(x_1))}(Q).
$$
\end{proposizione}
Let $(Q,\sigma)$ be a symmetric quiver and $(x,\sigma(x))$ a
sink-source admissible pair. For every $V\in Rep(Q)$, we define
the reflection functors
$$
C^+_{(x,\sigma(x))}V:=C^-_{\sigma(x)}C^+_x V\quad\textrm{and}\quad
C^-_{(x,\sigma(x))}V:=C^-_xC^+_{\sigma(x)} V.
$$
\begin{proposizione}\label{C+nabla}
Let $(Q,\sigma)$ be a symmetric quiver and $V\in Rep(Q)$. If
$(x,\sigma(x))$ is a sink-source admissible pair, then
$$
\nabla C^+_{(x,\sigma(x))}V=C^+_{(x,\sigma(x))}\nabla
V\quad\textrm{and}\quad\nabla
C^-_{(x,\sigma(x))}V=C^-_{(x,\sigma(x))}\nabla V.
$$
In particular for every $x$ admissible sink and $y$ admissible
source we have
$$
V=\nabla V \Leftrightarrow C^+_{(x,\sigma(x))}V=\nabla
C^+_{(x,\sigma(x))}V\Leftrightarrow C^-_{(x,\sigma(x))}V=\nabla
C^-_{(x,\sigma(x))}V.
$$
\end{proposizione}
\begin{proof}
It follows from definition of reflection functors for symmetric
quivers and from definition of duality functor $\nabla$.
\end{proof}
\begin{cor}\label{C+taunabla}
Let $(Q,\sigma)$ and $(Q',\sigma)$ be two symmetric quivers with
the same underlying graph. We suppose that
$Q'=c_{(x_m,\sigma(x_m))}\cdots c_{(x_1,\sigma(x_1))}(Q)$ for some
admissible sequence of sinks $x_1,\ldots,x_m\in Q_0$ for
admissible sink-source pairs and let
$V'=C^+_{(x_m,\sigma(x_m))}\cdots C^+_{(x_1,\sigma(x_1))}V\in
Rep(Q')$. Then
$$
V=C^-\nabla V\Leftrightarrow V'=C^-\nabla V'.
$$
\end{cor}
\begin{proof}
It follows from Proposition \ref{C+nabla}.
\end{proof}
\begin{proposizione}\label{sptosp}
Let $(Q,\sigma)$ be a symmetric quiver and let $x$ be an
admissible sink. Then $V$ is an orthogonal (resp. symplectic)
representation of $(Q,\sigma)$ if and only if
$C^+_{(x,\sigma(x))}V$ is an orthogonal (resp. symplectic) representation of $(c_{(x,\sigma(x))}Q,\sigma)$.\\
Similarly
for $C^-_{(x,\sigma(x))}$ if $x$ is an admissible source.
\end{proposizione}
\begin{proof} By proposition \ref{C+nabla} we have $V=\nabla V$
if and only if $C^+_{(x,\sigma(x))}V=\nabla C^+_{(x,\sigma(x))}V$.
To define an orthogonal (respectively symplectic) structure on
$C^+_{(x,\sigma(x))}V$ the only problem could occur at the
vertices fixed by $\sigma$. But, by definition of admissible sink
and of the involution $\sigma$, fixed vertices and fixed arrows
don't change under our reflection. The proof is similar for
$C^-_{(x,\sigma(x))}$ with $x$ an admissible source.
\end{proof}

\subsection{Orthogonal and symplectic semi-invariants}
If $W$ is a vector space of dimension $n$, we denote
$\widetilde{Gr}(r,W)$ the set of all decomposable tensors
$w_1\wedge\ldots\wedge w_r$, with $w_1,\ldots,w_r\in W$, inside
$\bigwedge^r W$.
\begin{lemma}\label{kac}
If $x$ is an admissible sink or source for a symmetric quiver
$(Q,\sigma)$ and $\alpha$ is a dimension vector such that
$c_{(x,\sigma(x))}\alpha(x)\geq 0$, then
\begin{itemize}
\item[i)] if $c_{(x,\sigma(x))}\alpha(x)> 0$ there exist isomorphisms
$$
SpSI(Q,\alpha)\stackrel{\varphi^{Sp}_{x,\alpha}}{\longrightarrow}
SpSI(c_{(x,\sigma(x))}Q,c_{(x,\sigma(x))}\alpha)
$$
and
$$
OSI(Q,\alpha)\stackrel{\varphi^{O}_{x,\alpha}}{\longrightarrow}
OSI(c_{(x,\sigma(x))}Q,c_{(x,\sigma(x))}\alpha),
$$
\item[ii)] if $c_{(x,\sigma(x))}\alpha(x)= 0$ there exist isomorphisms
$$
SpSI(Q,\alpha)\stackrel{\varphi^{Sp}_{x,\alpha}}{\longrightarrow}
SpSI(c_{(x,\sigma(x))}Q,c_{(x,\sigma(x))}\alpha)[y]
$$
and
$$
OSI(Q,\alpha)\stackrel{\varphi^{O}_{x,\alpha}}{\longrightarrow}
OSI(c_{(x,\sigma(x))}Q,c_{(x,\sigma(x))}\alpha)[y]
$$
\end{itemize}
where $A[y]$ denotes a polynomial ring with coefficients in $A$.
\end{lemma}
\begin{proof} We will prove the lemma for the symplectic case because the
orthogonal case is similar. Let $x\in Q_0$ be an admissible sink.
Put $r=\alpha(x)$ and $n=\sum_{ha=x}\alpha(ta)$. We note that
$c_{(x,\sigma(x))}\alpha(x)=n-r$. Put $V=\Bbbk^r$,
$V'=\Bbbk^{n-r}$ and
$W=\bigoplus_{ha=x}\Bbbk^{\alpha(ta)}\cong\Bbbk^n$. We define
$$Z=\bigoplus_{{a\in Q_1^+\atop ha\neq
x}}Hom(\Bbbk^{\alpha(ta)},\Bbbk^{\alpha(ha)})\oplus\bigoplus_{a\in
Q_1^{\sigma}}S^2(\Bbbk^{\alpha(ta)})^*$$ and
$$G=\prod_{{y\in
Q_0^+\atop y\neq x}}SL(\alpha(y))\times\prod_{y\in
Q_0^{\sigma}}Sp(\alpha(y)).
$$
\textit{Proof of i).} If $c_{(x,\sigma(x))}\alpha(x)> 0$ we have
$$
SpSI(Q,\alpha)=\Bbbk[SpRep(Q,\alpha)]^{SSp(Q,\alpha)}=
$$
$$
\Bbbk[Z\times Hom(W,V)]^{G\times SL\,V}=(\Bbbk[Z]\otimes\Bbbk[
Hom(W,V)]^{SL\,V}]^G=
$$
$$
(\Bbbk[Z]\otimes\Bbbk[\widetilde{Gr}(r,W^*)])^G
$$
and
$$
SpSI(c_{(x,\sigma(x))}Q,c_{(x,\sigma(x))}\alpha)=
$$
$$
\Bbbk[SpRep(c_{(x,\sigma(x))}Q,c_{(x,\sigma(x))}\alpha)]^{SSp(c_{(x,\sigma(x))}Q,c_{(x,\sigma(x))}\alpha)}=
$$
$$
\Bbbk[Z\times Hom(V',W)]^{G\times SL\,V'}=(\Bbbk[Z]\otimes\Bbbk[
Hom(V',W)]^{SL\,V'}]^G=
$$
$$
(\Bbbk[Z]\otimes\Bbbk[\widetilde{Gr}(n-r,W)])^G.
$$
Since $\widetilde{Gr}(r,W^*)$ and $\widetilde{Gr}(n-r,W)$ are
isomorphic as $G$-varieties, it follows that $SpSI(Q,\alpha)$ and
$SpSI(c_{(x,\sigma(x))}Q,c_{(x,\sigma(x))}\alpha)$ are
isomorphic.\\
\textit{Proof ii).} If $c_{(x,\sigma(x))}\alpha(x)= 0$, then $n=r$
and $V'=0$. So $\widetilde{Gr}(0,W)$ is a point and hence
\begin{equation}\label{era}
SpSI(Q,\alpha)=(\Bbbk[Z]\otimes \Bbbk[Hom(W,V)])^{G\times SL\,V}
\end{equation}
is isomorphic to
$$
SpSI(c_{(x,\sigma(x))}Q,c_{(x,\sigma(x))}\alpha)=(\Bbbk[Z]\otimes
\Bbbk[Hom(V',W)])^{G\times SL\,V}=\Bbbk[Z]^{G\times SL(V)}.
$$
Now let $A=\{a\in Q^+_1|\,ha=x\}$. Using theorem \ref{fc}, each
summand of (\ref{era}) contains\\ $(\bigotimes_{a\in
A}S_{\lambda(a)}V)^{SL\,V}$ as factor. By proposition \ref{i2}
each $\lambda(a)$, with $a\in A$, has to contain a column of
height $\alpha(ta)$, hence $\lambda(a)=\mu(a)+(1^{\alpha(ta)})$,
for some $\mu(a)$ in the set of partitions $\Lambda$. So as factor
we have
$$
\bigotimes_{a\in
A}(S_{(1^{\alpha(ta)})}\Bbbk^{\alpha(ta)})^{SL\,V_{ta}}\otimes\left(\bigotimes_{a\in
A}S_{(1^{\alpha(ta)})}V\right)^{SL\,V}
$$
which is generated by
$det(\bigoplus_{ha=x}\Bbbk^{\alpha(ta)}\rightarrow\Bbbk^{\alpha(x)})$.
On the other hand we have\\ $\Bbbk[Hom(W,V)])^{G\times
SL\,V}=\Bbbk[det(\bigoplus_{ha=x}\Bbbk^{\alpha(ta)}\rightarrow\Bbbk^{\alpha(x)})]$
and so we have the statement \textit{ii)}, with
$y=det(\bigoplus_{ha=x}\Bbbk^{\alpha(ta)}\rightarrow\Bbbk^{\alpha(x)})$.
\end{proof}
  We recall that, by definition, symplectic groups or orthogonal
groups act
  on the spaces which are defined on the
  vertices in $Q_0^{\sigma}$, so we have
  \begin{definizione}\label{wsQ}
  Let $V$ be a
  representation of the underlying quiver $Q$ with $\underline{dim}V=\alpha$ such that
  $\langle\alpha,\beta\rangle=0$ for some symmetric
  dimension vector $\beta$. The weight of $c^V$ on
  $SpRep(Q,\beta)$ (respectively on $ORep(Q,\beta)$) is
  $\langle\alpha,\cdot\rangle-\sum_{x\in
  Q_0^{\sigma}}\varepsilon_{x,\alpha}$, where
  \begin{equation}\label{epsilon}
  \varepsilon_{x,\alpha}(y)=\left\{\begin{array}{ll}
  \langle\alpha,\cdot\rangle(x) & y=x\\
  0 & \textrm{otherwise}.\end{array}\right.
  \end{equation}
  \end{definizione}
  We shall say that a weight is \textit{symmetric} if
  $\chi(i)=-\chi(\sigma(i))$ for every $i\in Q_0$.
  \begin{oss}\label{pesi}
  Let $(Q,\sigma)$ be a symmetric quiver and $V\in Rep(Q,\alpha)$. We note that
  $$
  \langle\underline{dim}(C^-\nabla
  V),\cdot\rangle(i)=-\langle\alpha,\cdot\rangle(\sigma(i))
  $$
  for every $i\in Q_0$. So, if $C^-\nabla V=V$ then
  $\chi=\langle\alpha,\cdot\rangle$ is a symmetric weight.
  \end{oss}
We show the relation between $c^V$ and $C^+_x$
(respectively $C^-_x$) and between $c^V$ and duality functor $\nabla$.\\
If $f$ is a semi-invariant of a quiver $Q$, we call $Z(f)$ the
 vanishing set of $f$. We have following lemma
 \begin{lemma}\label{zeridif}
Let $f$ and $f'$ be two semi-invariants of a quiver $Q$ such that
$Z(f)=Z(f')$ is
 irreducible. Then $f=k\cdot f'$ for some non zero
 $k\in\Bbbk$.
 \end{lemma}
\begin{proof} Since $Z(f)$ is irreducible, also $f$ is an
 irreducible polynomial. From $Z(f)=Z(f')$ it follows that $f'|f$
 and so $f=k\cdot f'$ for some non zero
 $k\in\Bbbk$. \end{proof}
 \begin{oss}\label{genirr}
 Let $\alpha$ be a dimension vector. For any finite set $S$ of generators of $SI(Q,\alpha)$, we can find a
finite generating set $S'$ consisting of irreducible polynomials
such that every element in $S$ is the product of elements in $S'$.
 \end{oss}
\begin{lemma}\label{cVW=c+xVW}
Let $V$ be an indecomposable representation of $Q$ of dimension
$\alpha$ such that $Z(c^V)$ is irreducible and let $x$ be a sink
of $Q$. Then
$$
c^V=k\cdot (c^{C^+_xV}\circ C^+_x)
$$
on $Rep(Q,\beta)$ such that $\langle\alpha,\beta\rangle=0$ and for
some non zero $k\in\Bbbk$.
\end{lemma}
\begin{proof} First we note that, by remark \ref{genirr} and by theorem
\ref{dw}, it's not restrictive to suppose $Z(c^V)$ is irreducible.
By remark \ref{qp}(ii), the vanishing set of $c^V$ is the
hypersurface
$$
Z(c^V)=\{W\in Rep(Q,\beta)|\, Hom_Q(V,W)\neq 0\}
$$
and the vanishing set of $c^{C^+_xV}$ is the hypersurface
$$
Z(c^{C^+_xV})=\{C^+_xW\in Rep(c_x(Q),c_x(\beta))|\,
Hom_Q(C^+_xV,C^+_xW)\neq 0\}.
$$
By definition of reflection functor, for every $W\in
Rep(Q,\beta)$,
$$
Hom_Q(V,W)\neq 0\Leftrightarrow Hom_Q(C^+_xV,C^+_xW)\neq 0.
$$
Hence $Z(c^V)=Z(c^{C^+_xV})$.\\
So, by lemma \ref{zeridif}, we conclude that there exist non zero
$k\in\Bbbk$ such that $c^V=k\cdot (c^{C^+_xV}\circ
C^+_x)$.\end{proof}

Similarly one proves the following two lemmas
\begin{lemma}\label{cVW=c-xVW}
Let $V$ be an indecomposable representation of $Q$ of dimension
$\alpha$ such that $Z(c^V)$ is irreducible and let $x$ be a source
of $Q$. Then
$$
c^V=k\cdot (c^{C^-_xV}\circ C^-_x)
$$
on $Rep(Q,\beta)$ such that $\langle\alpha,\beta\rangle=0$ and for
some non zero $k\in\Bbbk$.\hspace{5.4cm}$\square$
\end{lemma}
\begin{lemma}\label{cV=cVnabla1}
  Let $(Q,\sigma)$ be a symmetric quiver. For every representation
  $V$ of the underlying quiver $Q$ such that $Z(c^V)$ is irreducible, we have
  \begin{eqnarray}
  c^{V}=k\circ (c^{C^-\nabla V}\circ\nabla)
  \end{eqnarray}
  for some non zero $k\in\Bbbk$.\hspace{11.4cm}$\square$
  \end{lemma}
\begin{cor}\label{cV=cVnabla}
Let $(Q,\sigma)$ be a symmetric quiver. For every representation
  $V$ of the underlying quiver $Q$ and for every orthogonal or
  symplectic representation $W$ such that
  $\langle\underline{dim}(V),\underline{dim}(W)\rangle=0$, we have
  $$
  c^{V}(W)=c^{C^-\nabla V}(W).
  $$
\end{cor}
\begin{proof} It follows directly from lemma \ref{cV=cVnabla1}.
\end{proof}
We conclude this section with a lemma which will be useful later.
\begin{lemma}\label{cl}
Let
$$
\xymatrix@-1pc{&\ar@{.}[dr]&&&&\ar@{.}[dr]
&&&&\\
(Q,\sigma):&\ar@{.}[r]&y\ar[r]^{a}&x\ar[r]^{b}&z\ar@{.}[r]\ar@{.}[ur]\ar@{.}[dr]&
\ar@{.}[r]&\sigma(z)\ar[r]^{a}&\sigma(x)\ar[r]^{b}&\sigma(y)\ar@{.}[r]\ar@{.}[ur]\ar@{.}[dr]&\\
&\ar@{.}[ur]&&&&\ar@{.}[ur] &&&&}
$$
be a symmetric quiver. Assume there exist only two arrows in
$Q_1^+$ incident to $x\in Q_0^+$, $a:y\rightarrow x$ and
$b:x\rightarrow z$ with $y,z\in Q_0^+\cup Q_0^{\sigma}$. Let $V$
be an orthogonal or symplectic representation with symmetric dimension vector $(\alpha_i)_{i\in Q_0}=\alpha$ such that $\alpha_x\geq max\{\alpha_y,\alpha_z\}$.\\
We define the symmetric quiver $Q'=((Q_0',Q_1'),\sigma)$ with
$n-2$ vertices such that $Q_0'=Q_0\setminus\{x,\sigma(x)\}$ and
$Q_1'=Q_1\setminus\{a,b,\sigma(a),\sigma(b)\}\cup\{ba,\sigma(a)\sigma(b)\}$
and let $\alpha'$ be the dimension of $V$ restricted to $Q'$. \\
We have:
\begin{itemize}
\item[(Sp)] Assume $V$ symplectic. Then
\begin{itemize}
\item[(a)] if $\alpha_x> max\{\alpha_y,\alpha_z\}$ then
$SpSI(Q,\alpha)=SpSI(Q',\alpha')$,
\item[(b)] if $\alpha_x=\alpha_y>\alpha_z$ then
$SpSI(Q,\alpha)=SpSI(Q',\alpha')[detV(a)]$,
\item[(b')] if $\alpha_x=\alpha_z>\alpha_y$ then
$SpSI(Q,\alpha)=SpSI(Q',\alpha')[detV(b)]$,
\item[(c)] if $\alpha_x=\alpha_y=\alpha_z$ then
 $SpSI(Q,\alpha)=SpSI(Q',\alpha')[detV(a),detV(b)]$.
\end{itemize}
\item[(O)] Assume $V$ orthogonal. Then
\begin{itemize}
\item[(a)] if $\alpha_x> max\{\alpha_y,\alpha_z\}$ then
$OSI(Q,\alpha)=OSI(Q',\alpha')$,
\item[(b)] if $\alpha_x=\alpha_y>\alpha_z$ then
$OSI(Q,\alpha)=OSI(Q',\alpha')[detV(a)]$,
\item[(b')] if $\alpha_x=\alpha_z>\alpha_y$ then
$OSI(Q,\alpha)=OSI(Q',\alpha')[detV(b)]$,
\item[(c)] if $\alpha_x=\alpha_y=\alpha_z$ then
 $OSI(Q,\alpha)=OSI(Q',\alpha')[detV(a),detV(b)]$.
\end{itemize}
\end{itemize}
\end{lemma}
\begin{proof} By Cauchy formulas (Theorem \ref{fc}), we can
decompose $SpSI(Q,\alpha)$ and $OSI(Q,\alpha)$ in tensor products
of Schur modules. Now it's enough to adjust the proof of Lemma 28
in \cite{sw} to symmetric quivers using properties of Schur
modules (Proposition \ref{pms}) and Propositions \ref{i1} and
\ref{i2}.\end{proof}
\section{Semi-invariants of symmetric quivers of finite type}
We recall that by Theorem \ref{ctf} a symmetric quiver of finite
type has the underlying quiver of type $A_n$. Throughout this
section we enumerate vertices with $1,\ldots,n$ from left to right
and we call $a_i$ the arrow connecting $i$ and $i+1$, moreover we
define $\sigma$ by $\sigma(i)=n-i+1$, for every
$i\in\{1,\ldots,n\}$, and $\sigma(a_i)=\sigma(a_{n-i})$, for every
$i\in\{1,\ldots,n-1\}$.\\
We shall denote $V_{i,j}$ the indecomposable representations of
$A_n$ such that
$$
(\underline{dim}V_{i,j})_k=\left\{\begin{array}{lr}1&1\leq i\leq
k\leq j\leq n\\0&\textrm{otherwise}.\end{array}\right.
$$
Let $(Q,\sigma)$ be a symmetric quiver of type $A_n$ and let $V\in
Rep(Q)$ be indecomposable such that $C^-\nabla V=V$. By
Auslander-Reiten quiver of $Q$, we note that $V=V_{i,\sigma(i+1)}$
and the middle term $Z$ of the almost split sequence of $V$ is
$V_{i,\sigma(i)}\oplus V_{i+1,\sigma(i+1)}$, for some $i\in Q_0$.
Finally we note that we can define on $Z$ an orthogonal structure
if $n$ is odd and a symplectic structure if $n$ is even. So, for
an indecomposable representation of $Q$, if $n$ is even then
\textit{(i)} and \textit{(ii)} of property \textit{(Spp)} are
equivalent and if $n$ is odd then \textit{(i)} and \textit{(ii)}
of property \textit{(Op)} are
equivalent.\\
Let $(Q,\sigma)$ be a symmetric quiver and $V\in Rep(Q)$. We
consider
$$
0\longrightarrow
P_1\stackrel{d^V_{min}}{\longrightarrow}P_0\longrightarrow
V\longrightarrow 0
$$
the minimal projective resolution of $V$.
\begin{lemma}\label{ss2}
Let $(A_n,\sigma)$ be a symmetric quiver of type $A$. Let $V\in
Rep(Q,\alpha)$ such that $V=C^-\nabla V$ and let $\beta$ a
symmetric dimension vector such that
$\langle\alpha,\beta\rangle$=0, then we have the following.
\begin{itemize}
\item[(i)] If $n$ is even, $Hom_Q(d^V_{min},\cdot)$ is skew-symmetric on $ORep(Q,\beta)$.
\item[(ii)] If $n$ is odd, $Hom_Q(d^V_{min},\cdot)$ is skew-symmetric on $SpRep(Q,\beta)$.
\end{itemize}
\end{lemma}
\begin{proof}
Let $V\in Rep(Q,\alpha)$ be such that $C^-\nabla V=V$. So the
weight $\chi=\langle\alpha,\cdot\rangle-\sum_{x\in
Q_0^{\sigma}}\varepsilon_{x,\alpha}$ associated to $\alpha$ is
symmetric  (see Remark \ref{pesi}). If $m_1$ is the first vertex
such that $\chi(m_1)\neq 0$, in particular we suppose
$\chi(m_1)=1$, then the last vertex $m_s$ such that $\chi(m_s)\neq
0$ is $m_s=\sigma(m_1)$ and $\chi(m_s)=-1$. Between $m_1$ and
$m_s$, -1 and 1 alternate in correspondence respectively of sinks
and of sources. Moreover, by definition of symmetric weight, we
have $s=2l$ for some $l\in\mathbb{N}$. We call $i_2,\ldots,i_l$
the sources, $j_1,\ldots,j_{l-1}$ the sinks, $i_1=m_1$ and
$j_l=m_s$. Hence we have $\sigma(i_t)=j_{l-t+1}$ and
$i_1<j_1<\ldots<i_l<j_l$. Now the minimal projective resolution
for V is
\begin{equation}
0\longrightarrow\bigoplus_{j=j_1}^{j_l}
P_{j}\stackrel{d_{min}^V}{\longrightarrow}\bigoplus_{i=i_1}^{i_l}P_{i}\longrightarrow
V\longrightarrow 0
\end{equation}
and so we have
\begin{equation}
0\longrightarrow\bigoplus_{j=j_1}^{j_l}
P_{j}\stackrel{d_{min}^V}{\longrightarrow}\bigoplus_{j=j_1}^{j_l}P_{\sigma(j)}\longrightarrow
V\longrightarrow 0,
\end{equation}
with
\begin{equation}
(d_{min}^V)_{hk}=\left\{\begin{array}{ll} -a_{i_{k+1},j_k} &
\textrm{if}\quad
h=l-k\\
a_{i_k,j_k} & \textrm{if}\quad h=l-k+1\\
0 & \textrm{otherwise},
\end{array}\right.
\end{equation}
where $a_{i,j}$ is the oriented path from $i$ to $j$.\\
Hence
\begin{equation}
Hom(d_{min}^V,W):\bigoplus_{j=j_1}^{j_l}W(\sigma(j))=\bigoplus_{j=j_1}^{j_l}W(j)^*\longrightarrow\bigoplus_{j=j_1}^{j_l}W(j)
\end{equation}
where
\begin{equation}
(Hom(d_{min}^V,W))_{hk}=\left\{\begin{array}{ll}
-W(a_{i_{h+1},j_{h}}) &
\textrm{if}\quad k=l-h\\
W(a_{i_h,j_h}) &
\textrm{if}\quad k=l-h+1\\
0 & \textrm{otherwise}.
\end{array}\right.
\end{equation}
Now $W$ is orthogonal or symplectic, so for $k\neq h$, if
$k=l-h+1$ we have
$$
(Hom(d_{min}^V,W))_{hk}=W(a_{i_h,j_{h}})=W(a_{\sigma(j_{l-h+1}),j_{h}})=
-W(a_{\sigma(j_{h}),j_{l-h+1}})^t=
$$
$$
-W(a_{i_{l-h+1},j_{l-h+1}})^t=
-W(a_{i_k,j_k})^t=-((Hom(d_{min}^V,W))_{kh})^t.
$$
In a similar way it proves that if $k=l-h$ then
$(Hom(d_{min}^V,W))_{hk}=-((Hom(d_{min}^V,W))_{kh})^t$.\\
Finally the only cases for which $(Hom(d_{min}^V,W))_{hh}\neq 0$
are when $h=l-h+1$ and $h=l-h$. In the first case (the second one
is similar) we have
$(Hom(d_{min}^V,W))_{hh}=W(a_{i_h,j_{h}})=W(a_{\sigma(j_{h}),j_{h}})$
and $-((Hom(d_{min}^V,W))_{hh})^t=-W(a_{i_{h},j_{h}})^t=
-W(a_{\sigma(j_{h}),j_{h}})^t$. But
$W(a_{\sigma(j_{h}),j_{h}})=-W(a_{\sigma(j_{h}),j_{h}})^t$ for $n$
even if and only if $W\in ORep(Q)$, for $n$ odd if and only if
$W\in SpRep(Q)$.
\end{proof}
\begin{definizione}
For $V\in Rep(Q,\alpha)$ satisfying property \textit{(Spp)} (resp.
satisfying property \textit{(Op)} such that
$\langle\alpha,\beta\rangle=0$, where $\beta$ is an orthogonal
(resp. symplectic) dimension vector, we define
$$
\begin{array}{rcl}
pf^V:Rep(Q,\beta)&\longrightarrow&\Bbbk\\
W&\longmapsto&pf^V(W)=Pf(Hom_Q(d^V_{min},W)).
\end{array}
$$
\end{definizione}
By following Propositions and by Lemma \ref{kac}, it follows that
if Theorem \ref{tp1} and Theorem \ref{tp2} are true for
$(Q,\sigma)$, then they are true for
$(c_{(x,\sigma(x))}Q,\sigma)$.
\begin{proposizione}\label{SIcx1}
Let $(Q,\sigma)$ be a symmetric quiver of finite type. Let
$\alpha$ be a symmetric dimension vector, $x$ be an admissible
sink and $\varphi^{Sp}_{x,\alpha}$ be as defined in lemma \ref{kac}.\\
Then $\varphi^{Sp}_{x,\alpha}(c^V)=c^{C^+_{(x,\sigma(x))}V}$ and
$\varphi^{Sp}_{x,\alpha}(pf^W)=pf^{C^+_{(x,\sigma(x))}W}$, where
$V$ and $W$ are indecomposables of $Q$ such that $\langle
\underline{dim}\,V,\alpha
\rangle=0=\langle\underline{dim}\,W,\alpha \rangle$ and $W$
satisfies property \textit{(Op)}. In particular
\begin{itemize}
\item[(i)] if $0=\alpha_x\neq\alpha_{x-1}+\alpha_{x+1}$, then
$(\varphi^{Sp}_{x,\alpha})^{-1}(c^{S_x})=0$;
\item[(ii)] if $0\neq\alpha_x=\alpha_{x-1}+\alpha_{x+1}$, then
$\varphi^{Sp}_{x,\alpha}(c^{S_{\sigma(x)}})=0$.
\end{itemize}
\end{proposizione}
\begin{proof} We consider the same notation of proof of Lemma \ref{kac}.
If $x$ is an admissible sink of $(Q,\sigma)$, then we have, by
definition of $C^-_{(x,\sigma(x))}$,
$$C_{(x,\sigma(x))}^-|_{Z}(SpRep(c_{(x,\sigma(x))}Q,c_{(x,\sigma(x))}\alpha))=Z$$
and
$$C_{(x,\sigma(x))}^-|_{Hom(V',W)}(SpRep(c_{(x,\sigma(x))}Q,c_{(x,\sigma(x))}\alpha))=Hom(W,V).$$
Now $C_{(x,\sigma(x))}^-$ induces a ring morphism
$$
\begin{array}{rcl}
\phi_{x,\alpha}^{Sp}:\Bbbk[SpRep(Q,\alpha)]&\longrightarrow&\Bbbk[SpRep(c_{(x,\sigma(x))}Q,c_{(x,\sigma(x))}\alpha)]\\
f&\longmapsto& f\circ C_{(x,\sigma(x))}^-\end{array}
$$
By proof of Lemma \ref{kac}, we note that
$$
\Bbbk[C_{(x,\sigma(x))}^-Z\times
C_{(x,\sigma(x))}^-Hom(V',W)]^{SSp(Q,\alpha)}=\Bbbk[Z\times
Hom(W,V)]^{SSp(Q,\alpha)}
$$
is isomorphic by $\varphi_{x,\alpha}^{Sp}$ to $\Bbbk[Z\times
Hom(V',W)]^{SSp(c_{(x,\sigma(x))}Q,c_{(x,\sigma(x))}\alpha)}$.
Hence
$\varphi_{x,\alpha}^{Sp}=\phi_{x,\alpha}^{Sp}|_{SpSI(Q,\alpha)}$
and so for every representation $Z$ of dimension vector $\alpha$
of $(Q,\sigma)$ we have
\begin{equation}\label{varphi}
\varphi_{x,\alpha}^{Sp}(c^V)(C_{(x,\sigma(x))}^+Z)=(c^V\circ
C_{(x,\sigma(x))}^-)(C_{(x,\sigma(x))}^+Z)=c^V(Z)
\end{equation}
and
\begin{equation}\label{varphipf}\varphi_{x,\alpha}^{Sp}(pf^W)(C_{(x,\sigma(x))}^+Z)=(pf^W\circ
C_{(x,\sigma(x))}^-)(C_{(x,\sigma(x))}^+Z)=pf^W(Z).\end{equation}
By Lemma \ref{cVW=c+xVW} and \ref{cVW=c-xVW} we have
$c^V(Z)=k\cdot c^{C_{(x,\sigma(x))}^+V}(C_{(x,\sigma(x))}^+Z)$,
for some $k\in\Bbbk$. So, by (\ref{varphi}),
$\varphi_{x,\alpha}^{Sp}$ sends $c^V$ to
$c^{C_{(x,\sigma(x))}^+V}$ up to a constant in $\Bbbk$. Similarly
for $pf^W$. \end{proof}
\begin{proposizione}\label{SIcx2}
Let $(Q,\sigma)$ be a symmetric quiver of finite type. Let
$\alpha$ be a symmetric dimension vector, $x$ be an admissible
sink and $\varphi^{O}_{x,\alpha}$ be as defined in lemma \ref{kac}.\\
Then $\varphi^{O}_{x,\alpha}(c^V)=c^{C^+_{(x,\sigma(x))}V}$ and
$\varphi^{O}_{x,\alpha}(pf^W)=pf^{C^+_{(x,\sigma(x))}W}$, where
$V$ and $W$ are indecomposables of $Q$ such that $\langle
\underline{dim}\,V,\alpha
\rangle=0=\langle\underline{dim}\,W,\alpha \rangle$ and $W$
satisfies property \textit{(Spp)}. In particular
\begin{itemize}
\item[(i)] if $0=\alpha_x\neq\alpha_{x-1}+\alpha_{x+1}$, then
$(\varphi^{O}_{x,\alpha})^{-1}(c^{S_x})=0$;
\item[(ii)] if $0\neq\alpha_x=\alpha_{x-1}+\alpha_{x+1}$, then
$\varphi^{O}_{x,\alpha}(c^{S_{\sigma(x)}})=0$.
\end{itemize}
\end{proposizione}
\begin{proof} One proves similarly to Proposition
\ref{SIcx1}.
\end{proof}
In next sections we prove Theorem \ref{tp1} and \ref{tp2} for
symmetric quivers of type $A$ in equioriented case, i.e. the case
in which all the arrows have the same orientation. We shall call
$\Lambda$, $ER\Lambda$ and $EC\Lambda$ respectively the set of
partitions, the set of partition with even rows and the set of
partition of even columns.
\subsection{The symplectic case for $\boldsymbol{A_{2n}}$}
We can restate Theorem \ref{tp1} in the following way
\begin{teorema}\label{tfsij1}
Let $(Q,\sigma)$ be an equioriented symmetric quiver of type
$A_{2n}$ and let $\alpha$ be a symplectic dimension vector. Then
$SpSI(Q,\alpha)$ is generated by the following indecomposable
semi-invariants:
\begin{itemize}
\item[(i)] $c^{V_{j,i}}$ of weight $\langle\underline{dim}V_{j,i},\cdot\rangle$
for every $1\leq j\leq i\leq n-1$ such that
$\langle\underline{dim}\,V_{j,i},\alpha\rangle=0$,
\item[(ii)] $c^{V_{i,2n-i}}$ of weight
$\langle\underline{dim}V_{i,2n-i},\cdot\rangle$ for every
$i\in\{1,\ldots,n\}$.
\end{itemize}
\end{teorema}
The result follows from the following statement
\begin{teorema}\label{tfse}
Let $(Q,\sigma)$ be an equioriented symmetric quiver of type
$A_{2n}$, where
$$
Q=A_n^{eq}:1\stackrel{a_1}{\longrightarrow}2\cdots
n\stackrel{a_n}{\longrightarrow}n+1\cdots
2n-1\stackrel{a_{2n-1}}{\longrightarrow}2n,
$$
and let $V$ be a symplectic representation,
$\underline{dim}(V)=(\alpha_1,\ldots,\alpha_n)=\alpha$. Then
$SpSI(Q,\alpha)$ is generated by the following indecomposable
semi-invariants:
\begin{itemize}
\item[(i)] $det(V(a_i)\cdots V(a_j))$ with $j\leq i\in\{1,\ldots,n-1\}$ if $min\{\alpha_{j+1},\ldots,\alpha_{i}\}>\alpha_j=\alpha_{i+1}$;
\item[(ii)]  $det(V(a_{2n-i})\cdots V(a_i))$ with $i\in\{1,\ldots,n\}$ if $min\{\alpha_{i+1},\ldots,\alpha_{n}\}>\alpha_i.$
\end{itemize}
\end{teorema}
\begin{proof} We have
$$
 X:=SpRep(Q,\alpha)=\bigoplus_{i=1}^{n-1}V(ta_i)^*\otimes V(ha_i)\oplus S_2V_n^*.
$$
We proceed by induction on $n$.  For $n=1$ we have the symplectic
representation
$$
V_1\stackrel{V(a)}{\longrightarrow} V_1^*
$$
where $V_1$ is a vector space of dimension $\alpha$ and $V(a)$ is
a linear map such that $V(a)=V(a)^t$. So
$$
SpRep(Q,\alpha)=S^2V_1^*
$$
and by theorem \ref{fc}
$$
SpSI(Q,\alpha)=\bigoplus_{\lambda\in
ER\Lambda}(S_{\lambda}V_1)^{SL(V_1)},
$$
By proposition \ref{i1} and since $\lambda\in ER\Lambda$,
$SpSI(Q,\alpha)\neq 0$ if and only if
$\lambda=(\overbrace{2k,\ldots,2k}^{\alpha})$ for some
$k\in\mathbb{N}$ and we have that $(S_{\lambda}V_1)^{SL(V_1)}$ is
generated by a semi-invariant of weight $2k$. We note that
$V(a)\in S_2V_1^*\mapsto(detV(a))^{k}$ is a semi-invariant of
weight $2k$. So $SpSI(Q,\alpha)=\Bbbk[detV(a)]$.\\
 Now we prove the induction step. By theorem \ref{fc} we obtain
$$
SpSI(Q,\alpha)=\big(\Bbbk[X]\big)^{SL(V)}=
$$

$$
\bigoplus_{{\lambda(a_1),\ldots,\lambda(a_{n-1})\; and  \atop
\lambda(a_n)\in ER\Lambda}}(S_{\lambda(a_1)}V_1)^{SL(V_1)}\otimes
(S_{\lambda(a_1)}V_2^*\otimes
S_{\lambda(a_2)}V_2)^{SL(V_2)}\otimes \cdots\otimes
(S_{\lambda(a_{n-1})}V_n^*\otimes S_{\lambda(a_n)}V_n)^{SL(V_n)}.
$$
We suppose that there exists $i\in \{1,\dots,n-2\}$ such that
$\alpha_1\leq\cdots\leq\alpha_{i}$ and $\alpha_{i+i}<\alpha_i$. By
lemma \ref{cl},
$$
SpSI(Q,\alpha)=SpSI(Q^1,\alpha^1)
$$
where $Q^1$ is the smaller quiver $1\longrightarrow 2\cdots
i-1\longrightarrow i+1\cdots 2n-i+1\longrightarrow 2n-i+3\cdots
2n-1\longrightarrow 2n$ and $\alpha^1$ is the restriction of
$\alpha$ in $Q^1$.\\
If $i$ does't exist, we have $\alpha_1\leq\cdots\leq\alpha_{n-1}$.
So, by lemma \ref{cl}, we have the generators
$det\,V(a_i)=det\,V(\sigma(a_i))$ if
$\alpha_i=\alpha_{i+1}$, $1\leq i \leq n-2$.\\
We note that, by proposition \ref{i1},
$$
\lambda(a_1)=(\overbrace{k_1,\ldots,k_1}^{\alpha_1})
$$
is a rectangle with $k_1$ columns of height $\alpha_1$, for some
$k_1\in\mathbb{N}$. Since $\alpha_1\leq\cdots\leq\alpha_{n-1}$, by
proposition \ref{i2}, we obtain that there exist
$k_1,\ldots,k_{n-1}\in\mathbb{N}$ such that
$$
\lambda(a_{i})=(\overbrace{k_{i}+\cdots+k_1,\ldots,k_{i}+\cdots+k_1}^{\alpha_1},\ldots,\overbrace{k_{i},\ldots,k_{i}}^{\alpha_{i}-\alpha_{i-1}}),
$$
for every $i\in\{1,\ldots,n-1\}$.
 We also know that $\lambda_n$ must have even rows. If
 $\alpha_n=\alpha_j\leq \alpha_{j+1}\leq\cdots\leq \alpha_{n-1}$
 for some $j\in\{1,\ldots,n-1\}$ then $S_{\lambda_{n-1}}V_n^*=0$
 unless $ k_{n-1}+\cdots+k_{j+1}=0$, so
 $\lambda(a_{n-1})=\cdots=\lambda(a_{j+1})=\lambda(a_j)$. By proposition
 \ref{i2}, $(S_{\lambda(a_{n-1})}V_n^*\otimes S_{\lambda(a_n)}V_n)^{SL(V_n)}=(S_{\lambda(a_{j})}V_n^*\otimes
 S_{\lambda(a_n)}V_n)^{SL(V_n)}$ contains a semi-invariant if and
 only if
 $$
 \lambda(a_n)=(\overbrace{k_{n}+k_{j-1}+\cdots+k_1,\ldots,k_{n}+k_{j-1}+\cdots+k_1}^{\alpha_1},\ldots,\overbrace{k_{n},\ldots,k_{n}}^{\alpha_{n}-\alpha_{j-1}}),
 $$
 but
 $k_{n}+k_{j-1}+\cdots+k_1,k_{n}+k_{j-1}+\cdots+k_2,\ldots,k_n$
 have to be even and then $k_n,k_{j-1},\ldots,k_1$ have to be
 even. As before, by lemma \ref{cl}, we can consider the smaller
 quiver $Q^2:1\longrightarrow 2\cdots j\longrightarrow n\longrightarrow n+1 \longrightarrow 2n-j+1\cdots 2n-1\longrightarrow 2n$  and then
 $$
SpSI(Q,\alpha)\cong SpSI(Q^2,\alpha^2)=
$$
$$
(S_{\lambda(a_1)}V_1)^{SL(V_1)}\otimes\cdots\otimes(S_{\lambda(a_{j-1})}V_j^*\otimes
S_{\lambda(a_j)}V_j)^{SL(V_j)}\otimes(S_{\lambda(a_{j})}V_n^*\otimes
S_{\lambda(a_n)}V_n)^{SL(V_n)}.
$$
Now to complete the proof it's enough to find the generators of
$SpSI(Q^2,\alpha^2)$ for $\alpha_n=\alpha_j\leq
\alpha_{j+1}\leq\cdots\leq \alpha_{n-1}$.\\
By proposition \ref{i2}, for every $l\in
\{1,\ldots,j\},\;(S_{\lambda(a_{l-1})}V_l^*\otimes
S_{\lambda(a_l)}V_l)^{SL(V_l)}$ is generated by a semi-invariant
of weight $(0,\ldots,0,k_l,0,\dots,0)$ where $k_l=2h$ with
$h\in\mathbb{N}$, is $l$-th component. Since $V(a_{2n-l})\cdots
V(a_l)\in SpSI(Q,\alpha)\mapsto (det(V(a_{2n-l})\cdots V(a_l)))^h$
is a semi-invariant of weight $(0,\ldots,0,k_l,0,\dots,0)$, so it
generates $(S_{\lambda(a_{l-1})}V_l^*\otimes
S_{\lambda(a_l)}V_l)^{SL(V_l)}$. Now
$\lambda(a_l)=\lambda(a_{l-1})+(k_l^{\alpha_l})$ hence
$det(V(a_{2n-l})\cdots V(a_l))$ is a generator of
$SpSI(Q,\alpha)$.\\
In the summand of $SpSI(Q,\alpha)$ indexed by the families of
partitions in which
$\lambda(a_j)=(\overbrace{k_j,\ldots,k_j}^{\alpha_j=\alpha_n})$,
with $k_j\in\mathbb{N}$, we have that
$(S_{\lambda(a_{j})}V_j)^{SL(V_j)}\otimes
 (S_{\lambda(a_j)}V_n^*)^{SL(V_n)}$ is generated by
 a semi-invariant of weight\\
 $(0,\ldots,0,k_j,0,\dots,0,-k_j)$ where $k_j$ and $-k_j$ are
 respectively the $j$-th and the $n$-th component and we note, as
 before, that $(det(V(a_{n-1})\cdots V(a_j)))^{k_j}$ is a
 semi-invariant of weight $(0,\ldots,0,k_j,0,\dots,0,-k_j)$. Since
 $\lambda(a_j)=\lambda(a_{j-1})+(k_j^{\alpha_j=\alpha_n})$,
 $det(V(a_{n-1})\cdots V(a_j))$ is a generator of
 $SpSI(Q,\alpha)$.\\
In the summand of $SpSI(Q,\alpha)$ indexed by the families of
partitions in which
$\lambda(a_n)=(\overbrace{k_n,\ldots,k_n}^{\alpha_n})$ with
$k_n\in 2 \mathbb{N}$, we note again that
 $(S_{\lambda(a_n)}V_n)^{SL(V_n)}$ is generated by
 $(det(V(a_n)))^{k_n}$ of weight $(0,\ldots,0,k_n)$ where $n$-th component $k_n$ is
 even. Since $\lambda(a_n)=\lambda(a_{j-1})+(k_n^{\alpha_n})$, $det(V(a_n))$ is a generator of
 $SpSI(Q,\alpha)$.\end{proof}

\begin{proof}[Proof of Theorem \ref{tfsij1}] First we note that
$\alpha_j=\alpha_{i+1}$ is equivalent to
$\langle\underline{dim}\,V_{j,i},\underline{dim}\, V\rangle=0$. If
we consider the minimal projective resolution of $V_{j,i}$, we
have
$$
0\longrightarrow P_{i+1}\stackrel{a_i\cdots a_j}{\longrightarrow}
P_j\longrightarrow V_{j,i}\longrightarrow 0
$$
and applying the $Hom_Q$-functor we have
$$
Hom_Q(a_i\cdots a_j,V):Hom_Q(P_j,V)=V_j\stackrel{V(a_i\cdots
a_j)}{ \longrightarrow} V_{i+1}=Hom_Q(P_{i+1},V).
$$
So, $det(V(a_i)\cdots V(a_j))=det(Hom_Q(a_i\cdots
a_j,V))=c^{V_{j,i}}(V)$.\\
In the same way one proves that $det(V(a_{2n-i})\cdots
V(a_i))=det(V_i\longrightarrow
V_{2n-i+1}=V_i^*)=c^{V_{i,2n-i}}(V)$, but in this case, since
$\underline{dim}V=\underline{dim}\nabla V$, we have
$\alpha_i=\alpha_{2n-i+1}$ and so\\
$\langle\underline{dim}V_{i,2n-i},\underline{dim}V\rangle=0$ for
every $i\in\{1,\ldots,n\}$. Moreover we note that
\begin{itemize}
\item[(i)] $c^{V_{2n-i,2n-j}}(V)=c^{V_{j,i}}(V)$, by lemma \ref{cV=cVnabla}, since
$\tau^-\nabla V_{j,i}=V_{2n-i,2n-j}$;
\item[(ii)] for every $j\in\{1,\ldots,n-1\}$ and for every $i\in\{n+1,\ldots,2n-1\}\setminus\{2n-j\}$ there exists $k>j\in\{1,\ldots,n-1\}$ such that $2n-k=i$ and so
$c^{V_{j,i}}(V)=c^{V_{j,k-1}}(V)\cdot c^{V_{k,2n-k}}(V)$.
\end{itemize}
Now, using theorem \ref{tfse}, we obtain the statement of the
theorem. \end{proof}

\subsection{The orthogonal case for $\boldsymbol{A_{2n}}$}
We restate Theorem \ref{tp2} in the following way
\begin{teorema}\label{tfoij1}
Let $(Q,\sigma)$ be an equioriented symmetric quiver of type
$A_{2n}$ and let $\alpha$ be the dimension vector of an orthogonal representation of $(Q,\sigma)$.\\
Then $OSI(Q,\alpha)$ is generated by the following indecomposable
semi-invariants:
\begin{itemize}
\item[(i)] $c^{V_{j,i}}$ of weight $\langle\underline{dim}V_{j,i},\cdot\rangle$
for every $1\leq j\leq i\leq n-1$ such that
$\langle\underline{dim}\,V_{j,i},\alpha\rangle=0$,
\item[(ii)] $pf^{V_{i,2n-i}}$ of weight
$\frac{\langle\underline{dim}V_{i,2n-i},\cdot\rangle}{2}$ for
every $i\in\{1,\ldots,n\}$ such that $\alpha_i$ is even.
\end{itemize}
\end{teorema}
The result follows from the following statement
\begin{teorema}\label{tfoe}
Let $(Q,\sigma)$ be an equioriented symmetric quiver of type
$A_{2n}$, where
$$
Q=A_n^{eq}:1\stackrel{a_1}{\longrightarrow}2\cdots
n\stackrel{a_n}{\longrightarrow}n+1\cdots
2n-1\stackrel{a_{2n-1}}{\longrightarrow}2n,
$$
and let $V$ be an orthogonal representation,
$\underline{dim}(V)=(\alpha_1,\ldots,\alpha_n)=\alpha$. Then
$OSI(Q,\alpha)$ is generated by the following indecomposable
semi-invariants:
\begin{itemize}
\item[(i)] $det(V(a_i)\cdots V(a_j))$ with $j\leq i\in\{1,\ldots,n-1\}$ if $min(\alpha_{j+1},\ldots,\alpha_{i})>\alpha_j=\alpha_{i+1}$;
\item[(ii)]  $pf(V(a_{2n-i})\cdots V(a_i))$ with $i\in\{1,\ldots,n\}$ if
$min(\alpha_{i+1},\ldots,\alpha_{n})>\alpha_i$and $\alpha_i$ is
even.
\end{itemize}
\end{teorema}
\begin{proof} We have
$$
 X:=ORep(Q,\alpha)=\bigoplus_{i=1}^{n-1}V(ta_i)^*\otimes V(ha_i)\oplus \bigwedge^2V_n^*.
$$
We proceed by induction on $n$.  For $n=1$ we have the orthogonal
representation
$$
V_1\stackrel{V(a)}{\longrightarrow} V_1^*
$$
where $V_1$ is a vector space of dimension $\alpha$ and $V(a)$ is
a linear map such that $V(a)=-V(a)^t$.
$$
ORep(Q,\alpha)=\bigwedge^2V_1^*
$$
and by theorem \ref{fc}
$$
OSI(Q,\alpha)=\bigoplus_{\lambda\in
EC\Lambda}(S_{\lambda}V_1)^{SL(V_1)}
$$
By proposition \ref{i1} since $\lambda\in EC\Lambda$,
$OSI(Q,\alpha)\neq 0$ if and only if
$\lambda=(\overbrace{k,\ldots,k}^{\alpha})$ with $\alpha$ even,
for some $k$. We note that $V(a)\in\bigwedge^2V^*\mapsto
(pfV(a))^k$ is a semi-invariant of weight $k$ so
$(S_{\lambda}V_1)^{SL(V_1)}$
is generated by the semi-invariant $(pf\,V(a))^k$ if $\alpha$ is even and $OSI(Q,\alpha)=\Bbbk[pfV(a)]$.\\
Now we prove the induction step. By theorem \ref{fc} we obtain
$$
OSI(Q,\alpha)=\big(\Bbbk[X]\big)^{SL(V)}=
$$

$$
\bigoplus_{{\lambda(a_1),\ldots,\lambda(a_{n-1})\; and  \atop
\lambda(a_n)\in
EC\Lambda}}\!\!\!\!\!\!(S_{\lambda(a_1)}V_1)^{SL(V_1)}\otimes
(S_{\lambda(a_1)}V_2^*\otimes
S_{\lambda(a_2)}V_2)^{SL(V_2)}\otimes \cdots\otimes
(S_{\lambda(a_{n-1})}V_n^*\otimes S_{\lambda(a_n)}V_n)^{SL(V_n)}
$$
The proof of this theorem is the same of the proof of the theorem
\ref{tfse} up to when we have to consider $\alpha_n$. As in the
previous proof we can suppose $\alpha_1\leq\cdots\leq
\alpha_{n-1}$, otherwise, by induction, we can reduce to a smaller
quiver.\\
By lemma \ref{cl}, we have the generators
$det\,V(a_i)=det\,V(\sigma(a_i))$ if
$\alpha_i=\alpha_{i+1}$, $1\leq i\leq n-2$.\\
By proposition \ref{i2}, we obtain that there exist
$k_1,\ldots,k_{n-1}\in\mathbb{N}$ such that
$$
\lambda(a_{i})=(\overbrace{k_{i}+\cdots+k_1,\ldots,k_{i}+\cdots+k_1}^{\alpha_1},\ldots,\overbrace{k_{i},\ldots,k_{i}}^{\alpha_{i}-\alpha_{i-1}}),
$$
for every $i\in\{1,\ldots,n-1\}$.\\
Now we consider the hypothesis on $\lambda(a_n)$ by which it must
have even columns. If
 $\alpha_n=\alpha_j\leq \alpha_{j+1}\leq\cdots\leq \alpha_{n-1}$
 for some $j\in\{1,\ldots,n-1\}$ then $S_{\lambda(a_{n-1})}V_n^*=0$
 unless $ k_{n-1}+\cdots+k_{j+1}=0$, so
 $\lambda(a_{n-1})=\cdots=\lambda(a_{j+1})=\lambda(a_j)$. By proposition
 \ref{i2}, $(S_{\lambda(a_{n-1})}V_n^*\otimes S_{\lambda(a_n)}V_n)^{SL(V_n)}=(S_{\lambda(a_{j})}V_n^*\otimes
 S_{\lambda(a_n)}V_n)^{SL(V_n)}$ contains a semi-invariant if and
 only if
 $$
 \lambda(a_n)=(\overbrace{k_{n}+k_{j-1}+\cdots+k_1,\ldots,k_{n}+k_{j-1}+\cdots+k_1}^{\alpha_1},\ldots,\overbrace{k_{n},\ldots,k_{n}}^{\alpha_{n}-\alpha_{j-1}}),
 $$
 but
 $\alpha_1,\alpha_2-\alpha_1,\ldots,\alpha_n-\alpha_{j-1}$
 have to be even and then $\alpha_1,\ldots,\alpha_{j-1},\alpha_n$ have to be
 even. As before, by lemma \ref{cl}, we can consider the smaller
 quiver $Q^1:1\longrightarrow 2\cdots j\longrightarrow n\longrightarrow n+1 \longrightarrow 2n-j+1\cdots 2n-1\longrightarrow 2n$  and then
 $$
OSI(Q,\alpha)\cong
(S_{\lambda(a_1)}V_1)^{SL(V_1)}\otimes\cdots\otimes(S_{\lambda(a_{j-1})}V_j^*\otimes
S_{\lambda(a_j)}V_j)^{SL(V_j)}\otimes
(S_{\lambda(a_{j})}V_n^*\otimes S_{\lambda(a_n)}V_n)^{SL(V_n)}.
$$
Now to complete the proof it's enough to find the generator of
this algebra for $\alpha_n=\alpha_j\leq \alpha_{j+1}\leq\cdots\leq
\alpha_{n-1}$.\\
By proposition \ref{i2}, for every $l\in \{1,\ldots,j\}$ such that
$\alpha_l$ is even, $(S_{\lambda(a_{l-1})}V_l^*\otimes
S_{\lambda(a_l)}V_l)^{SL(V_l)}$ is generated by a semi-invariant
of weight $(0,\ldots,0,k_l,0,\dots,0)$ where $k_l\in\mathbb{N}$,
is $l$-th component. We note that $V(a_{2n-l})\cdots V(a_l)\in
OSI(Q,\alpha)\mapsto (pf(V(a_{2n-l})\cdots V(a_l)))^{k_l}$ is a
semi-invariant of weight
 $(0,\ldots,0,k_l,0,\dots,0)$, so it
generates\\ $(S_{\lambda(a_{l-1})}V_l^*\otimes
S_{\lambda(a_l)}V_l)^{SL(V_l)}$. Since
$\lambda(a_l)=\lambda(a_{l-1})+(k_l^{\alpha_l})$,\\
$pf(V(a_{2n-l})\cdots V(a_l))$ is a generator of $OSI(Q,\alpha)$.
In the summand of $OSI(Q,\alpha)$ indexed by the families of
partitions in which
$\lambda(a_j)=(\overbrace{k_j,\ldots,k_j}^{\alpha_j=\alpha_n})$,
with $k_j\in\mathbb{N}$, we have that
$(S_{\lambda(a_{j})}V_j)^{SL(V_j)}\otimes
 (S_{\lambda(a_j)}V_n^*)^{SL(V_n)}$ is generated by
 a semi-invariant of weight\\
 $(0,\ldots,0,k_j,0,\dots,0,-k_j)$ where $k_j$ and $-k_j$ are
 respectively the $j$-th and the $n$-th component and we note, as
 before, that $(det(V(a_{n-1})\cdots V(a_j)))^{k_j}$ is a
 semi-invariant of weight
 $(0,\ldots,0,k_j,0,\dots,0,-k_j)$. Since
 $\lambda(a_j)=\lambda(a_{j-1})+(k_j^{\alpha_j=\alpha_n})$,
 $det(V(a_{n-1})\cdots V(a_j))$ is a generator of
 $OSI(Q,\alpha)$.\\
In the summand of $OSI(Q,\alpha)$ indexed by the families of
partitions in which
$\lambda(a_n)=(\overbrace{k_n,\ldots,k_n}^{\alpha_n})$ with
$k_n\in \mathbb{N}$, we note again that if $\alpha_n$ is even
 $(S_{\lambda(a_n)}V_n)^{SL(V_n)}$ is generated by
 $(pf(V(a_n)))^{k_n}$ of weight $(0,\ldots,0,k_n)$. Since $\lambda(a_n)=\lambda(a_{j-1})+(k_n^{\alpha_n})$, $pf(V(a_n))$ is a generator of
 $SpSI(Q,\alpha)$. \end{proof}
\begin{proof}[Proof of Theorem \ref{tfoij1}] By lemma \ref{ss2}, we can
define $pf^V$ if $V=\tau^-\nabla V$, since we are dealing with
orthogonal case. Moreover we note that $V_{i, 2n-i}=\tau^-\nabla
V_{i, 2n-i}$. Hence using the theorem \ref{tfoe}, the  proof is
similar to the proof of theorem \ref{tfsij1}. \end{proof}

\subsection{The symplectic case for $\boldsymbol{A_{2n+1}}$}
We can restate theorem \ref{tp1} in the following way
\begin{teorema}\label{tfsij2}
Let $(Q,\sigma)$ be an equioriented symmetric quiver of type
$A_{2n+1}$ and let $\alpha$ be a symplectic dimension vector. Then
$SpSI(Q,\alpha)$ is generated by the following indecomposable
semi-invariants:
\begin{itemize}
\item[(i)] $c^{V_{j,i}}$ of weight
$\langle\underline{dim}V_{j,i},\cdot\rangle-\varepsilon_{n+1,\underline{dim}\,V_{j,i}}$
(see Definition \ref{wsQ}) for every $1\leq j\leq i\leq n$ such
that $\langle\underline{dim}\,V_{j,i},\alpha\rangle=0$.
\item[(ii)] $pf^{V_{i,2n+1-i}}$ of weight
$\frac{\langle\underline{dim}V_{i,2n+1-i},\cdot\rangle}{2}$ for
every $i\in\{1,\ldots,n\}$ such that $\alpha_i$ is even.
\end{itemize}
\end{teorema}
The result follows from the following statement
\begin{teorema}\label{tfse1}
Let $(Q,\sigma)$ be an equioriented symmetric quiver of type
$A_{2n+1}$, where
$$
Q:1\stackrel{a_1}{\longrightarrow}2\cdots
n\stackrel{a_n}{\longrightarrow}n+1\stackrel{a_{n+1}}{\longrightarrow}n+2\cdots
2n\stackrel{a_{2n}}{\longrightarrow}2n+1,
$$
and let $V$ be an symplectic representation,
$\underline{dim}(V)=(\alpha_1,\ldots,\alpha_{n+1})=\alpha$. Then
$SpSI(Q,\alpha)$ is generated by the following indecomposable
semi-invariants:
\begin{itemize}
\item[(i)] $det(V(a_i)\cdots V(a_j))$ with $j\leq i\in\{1,\ldots,n+1\}$ if $min(\alpha_{j+1},\ldots,\alpha_{i})>\alpha_j=\alpha_{i+1}$;
\item[(ii)]  $pf(V(a_{2n-i+1})\cdots V(a_i))$ with $i\in\{1,\ldots,n\}$ if
$min(\alpha_{i+1},\ldots,\alpha_{n+1})>\alpha_i$ and $\alpha_i$ is
even.
\end{itemize}
\end{teorema}
\begin{proof} First we recall that if $V$ is a symplectic representation
of dimension $\alpha=(\alpha_1,\ldots,\alpha_{n+1})$ of a
symmetric quiver of type $A_{2n+1}$, in the symplectic case,
  $V_{n+1}=V_{n+1}^*$ is a symplectic space, so if $V_{n+1}\neq 0$ then $dim\,V_{n+1}$ has to be
  even. We proceed by induction on $n$.  For $n=1$ we have
the symplectic representation
$$
V_1\stackrel{V(a)}{\longrightarrow}
V_2=V_2^*\stackrel{-V(a)^t}{\longrightarrow} V_1^*.
$$
By theorem \ref{fc}
$$
 SpSI(Q,\alpha)=\bigoplus_{\lambda\in
\Lambda}(S_{\lambda}V_1)^{SL(V_1)}\otimes(S_{\lambda}V_2)^{Sp(V_2)}.
$$
By proposition \ref{i1} and proposition \ref{i3},
$SpSI(Q,\alpha)\neq 0$ if and only if
$\lambda=(\overbrace{k,\ldots,k}^{\alpha_1})$, for some $k$, and
$ht(\lambda)$ has to be even. Moreover we have that
$(S_{(k^{\alpha_1})}V_1)^{SL(V_1)}\otimes(S_{(k^{\alpha_1})}V_2)^{Sp(V_2)}$
is generated by a semi-invariant of weight $(k,0)$. If
$\alpha_1>\alpha_2$ then $S_{\lambda}V_2=0$ unless $\lambda=0$ and
in this case $SpSI(Q,\alpha)=\Bbbk$. If $\alpha_1=\alpha_2$ then
$ht(\lambda)=\alpha_1=\alpha_2$. We note that $detV(a)^k$ is a
semi-invariant of weight $(k,0)$. Hence
$(S_{\lambda}V_1)^{SL(V_1)}\otimes(S_{\lambda}V_2)^{Sp(V_2)}$ is
generated by the semi-invariant $detV(a)^k$, so
$SpSI(Q,\alpha)=\Bbbk[detV(a)]$. Finally if $\alpha_1<\alpha_2$
then $ht(\lambda)=\alpha_1$ has to be even. We recall that in the
symplectic case $-V(a)^tV(a)$ is skew-symmetric. We note that
$pf(-V(a)^tV(a))^k$ is a semi-invariant of weight $(k,0)$ so
$(S_{\lambda}V_1)^{SL(V_1)}\otimes(S_{\lambda}V_2)^{Sp(V_2)}$ is
generated by the semi-invariant $pf(-V(a)^tV(a))^k$ if $\alpha_1$
is even and thus $SpSI(Q,\alpha)=\Bbbk[pf(-V(a)^tV(a))].$\\
 Now we prove the
induction step. Let $X=SpRep(Q,\alpha)$
 and by theorem \ref{fc} we obtain
$$
SpSI(Q,\alpha)=\big(\Bbbk[X]\big)^{SSp(V)}=
$$

$$
\bigoplus_{{\lambda(a_1),\ldots,\lambda(a_{n})\in
\Lambda}}(S_{\lambda(a_1)}V_1)^{SL(V_1)}\otimes
(S_{\lambda(a_1)}V_2^*\otimes
S_{\lambda(a_2)}V_2)^{SL(V_2)}\otimes
$$
$$
\cdots\otimes (S_{\lambda(a_{n-1})}V_n^*\otimes
S_{\lambda(a_n)}V_n)^{SL(V_n)}\otimes
(S_{\lambda(a_n)}V_{n+1})^{Sp(V_{n+1})} ,
$$
The proof of this theorem is the same of the proof of the theorem
\ref{tfse} up to when we have to consider $\alpha_{n+1}$. As in
the proof of theorem \ref{tfse} we can suppose
$\alpha_1\leq\cdots\leq \alpha_{n}$, otherwise, by induction, we
can reduce to a smaller quiver.\\
By lemma \ref{cl}, we have the generators
$det\,V(a_i)=det\,V(\sigma(a_i))$ if
$\alpha_i=\alpha_{i+1}$, $1\leq i\leq n-1$.\\
By proposition \ref{i2}, we obtain that there exist
$k_1,\ldots,k_{n}\in\mathbb{N}$ such that
$$
\lambda(a_{i})=(\overbrace{k_{i}+\cdots+k_1,\ldots,k_{i}+\cdots+k_1}^{\alpha_1},\ldots,\overbrace{k_{i},\ldots,k_{i}}^{\alpha_{i}-\alpha_{i-1}}),
$$
for every $i\in\{1,\ldots,n\}$.\\
Now, by proposition \ref{i3}, $\lambda(a_n)$ must have even
columns. If
 $\alpha_{n+1}=\alpha_j\leq \alpha_{j+1}\leq\cdots\leq \alpha_{n}$
 for some $j\in\{1,\ldots,n\}$ then $S_{\lambda(a_{n})}V_{n+1}^*=0$
 unless $ k_{n}+\cdots+k_{j+1}=0$, so
 $\lambda(a_{n})=\cdots=\lambda(a_{j+1})=\lambda(a_j)$. As before, by lemma \ref{cl}, we can consider the smaller
 quiver $Q^1:1\longrightarrow 2\cdots j\longrightarrow n+1\longrightarrow 2n-j+2\cdots 2n\longrightarrow 2n+1$  and then
 $$
SpSI(Q,\alpha)\cong
(S_{\lambda(a_1)}V_1)^{SL(V_1)}\otimes\cdots\otimes(S_{\lambda(a_{j-1})}V_j^*\otimes
S_{\lambda(a_j)}V_j)^{SL(V_j)}\otimes
$$
\begin{equation}\label{01}
(S_{\lambda(a_{j-1})}V_n^*\otimes
S_{\lambda(a_j)}V_n)^{SL(V_n)}\otimes(S_{\lambda(a_j)}V_{n+1})^{Sp(V_{n+1})},
\end{equation}
where
$$
\lambda(a_{j})=(\overbrace{k_{j}+\cdots+k_1,\ldots,k_{j}+\cdots+k_1}^{\alpha_1},\ldots,\overbrace{k_{j},\ldots,k_{j}}^{\alpha_{n+1}-\alpha_{j-1}}),
$$
and $\alpha_1,\alpha_2-\alpha_1,\dots,\alpha_{n+1}-\alpha_{j-1}$
have to be even otherwise, by proposition \ref{i3},
$(S_{\lambda(a_j)}V_{n+1})^{Sp(V_{n+1})}=0$.
 Now to complete the proof it's enough to find the generators of the algebra
 (\ref{01}) for $\alpha_{n+1}=\alpha_j\leq \alpha_{j+1}\leq\cdots\leq \alpha_{n}$.\\
By proposition \ref{i2}, for every $l\in \{1,\ldots,j\}$ such that
$\alpha_l$ is even, $(S_{\lambda(a_{l-1})}V_l^*\otimes
S_{\lambda(a_l)}V_l)^{SL(V_l)}$ is generated by a semi-invariant
of weight $(0,\ldots,0,k_l,0,\dots,0)$ where $k_l\in\mathbb{N}$,
is $l$-th component. We note that $V(a_{2n-l+1})\cdots V(a_l)\in
SpSI(Q,\alpha)\mapsto (pf(V(a_{2n-l+1})\cdots V(a_l)))^{k_l}$ is a
semi-invariant of weight
 $(0,\ldots,0,k_l,0,\dots,0)$, so it
generates $(S_{\lambda(a_{l-1})}V_l^*\otimes
S_{\lambda(a_l)}V_l)^{SL(V_l)}$. Since
$\lambda(a_l)=\lambda(a_{l-1})+(k_l)^{\alpha_l}$, then
$pf(V(a_{2n-l+1})\cdots V(a_l))$ is a generator of
$SpSI(Q,\alpha)$.\\
In the summand of $SpSI(Q,\alpha)$ indexed by the families of
partitions in which
$\lambda(a_j)=(\overbrace{k_j,\ldots,k_j}^{\alpha_j=\alpha_{n+1}})$,
with $k_j\in\mathbb{N}$, we have that
$(S_{\lambda(a_{j})}V_j)^{SL(V_j)}\otimes
 (S_{\lambda(a_j)}V_{n+1})^{Sp(V_{n+1})}$ is generated by
 a semi-invariant of weight
 $(0,\ldots,0,k_j,0,\dots,0,0)$ where $k_j$ is
  the $j$-th component and we note, as
 before, that $(det(V(a_{n})\cdots V(a_j)))^{k_j}$ is a
 semi-invariant of weight
 $(0,\ldots,0,k_j,0,\dots,0,0)$. Since
 $\lambda(a_j)=\lambda(a_{j-1})+(k_j)^{\alpha_j=\alpha_{n+1}}$,
 $det(V(a_{n})\cdots V(a_j))$ is a generator of $SpSI(Q,\alpha)$.
\end{proof}
\begin{proof}[Proof of Theorem \ref{tfsij2}] By lemma \ref{ss2}, we can
define $pf^V$ if $V=\tau^-\nabla V$, since we are dealing with
symplectic case. Moreover we note that $V_{i, 2n+1-i}=\tau^-\nabla
V_{i, 2n+1-i}$, for every $i\in\{1,\ldots,n\}$ . Hence using the
theorem \ref{tfse1}, the proof is similar to the proof of theorem
\ref{tfsij1}. \end{proof}

\subsection{The orthogonal case for $\boldsymbol{A_{2n+1}}$}
We restate the Theorem \ref{tp2} in the following way
\begin{teorema}\label{tfoij2}
Let $(Q,\sigma)$ be an equioriented symmetric quiver of type
$A_{2n+1}$ and let $\alpha$ be an orthogonal dimension vector.
Then $OSI(Q,\alpha)$ is generated by the following indecomposable
semi-invariants:
\begin{itemize}
\item[(i)] $c^{V_{j,i}}$ of weight
$\langle\underline{dim}V_{j,i},\cdot\rangle-\varepsilon_{n+1,\underline{dim}\,V_{j,i}}$
for every $1\leq j\leq i\leq n$ such that
$\langle\underline{dim}\,V_{j,i},\alpha\rangle=0$.
\item[(ii)] $c^{V_{i,2n+1-i}}$ of weight
$\langle\underline{dim}V_{i,2n+1-i},\cdot\rangle$ for every
$i\in\{1,\ldots,n\}$.
\end{itemize}
\end{teorema}
The result follows from the following statement
\begin{teorema}\label{tfoe1}
Let $(Q,\sigma)$ be an equioriented symmetric quiver of type
$A_{2n+1}$, where
$$
Q:1\stackrel{a_1}{\longrightarrow}2\cdots
n\stackrel{a_n}{\longrightarrow}n+1\stackrel{a_{n+1}}{\longrightarrow}n+2\cdots
2n\stackrel{a_{2n}}{\longrightarrow}2n+1,
$$
and let $V$ be an orthogonal representation,
$\underline{dim}(V)=(\alpha_1,\ldots,\alpha_{n+1})=\alpha$. Then
$OSI(Q,\alpha)$ is generated by the following indecomposable
semi-invariants:
\begin{itemize}
\item[(i)] $det(V(a_i)\cdots V(a_j))$ with $j\leq i\in\{1,\ldots,n+1\}$ if $min(\alpha_{j+1},\ldots,\alpha_{i})>\alpha_j=\alpha_{i+1}$;
\item[(ii)]  $det(V(a_{2n-i+1})\cdots V(a_i))$ with $i\in\{1,\ldots,n\}$ if
$min(\alpha_{i+1},\ldots,\alpha_{n+1})>\alpha_i$.
\end{itemize}
\end{teorema}
\begin{proof} We proceed by induction on $n$.  For $n=1$ we have the
orthogonal representation
$$
V_1\stackrel{V(a)}{\longrightarrow}
V_2=V_2^*\stackrel{-V(a)^t}{\longrightarrow} V_1^*
$$
where $V_1$ is a vector space of dimension $\alpha_1$, $V_2$ is a
orthogonal space of dimension $\alpha_2$ and $V(a)$ is a linear
map. By theorem \ref{fc}
$$
 OSI(Q,\alpha)=\bigoplus_{\lambda\in
\Lambda}(S_{\lambda}V_1)^{SL(V_1)}\otimes(S_{\lambda}V_2)^{SO(V_2)}.
$$
By proposition \ref{i1} and proposition \ref{i3},
$OSI(Q,\alpha)\neq 0$ if and only if
$\lambda=(\overbrace{k,\ldots,k}^{\alpha_1})$, for some $k\in
2\mathbb{N}$. Moreover we have that
$(S_{(k^{\alpha_1})}V_1)^{SL(V_1)}\otimes(S_{(k^{\alpha_1})}V_2)^{SO(V_2)}$
is generated by a semi-invariant of weight $(k,0)$. If
$\alpha_1>\alpha_2$ then $S_{\lambda}V_2=0$ unless $\lambda=0$ and
in this case $OSI(Q,\alpha)=\Bbbk$. If $\alpha_1=\alpha_2$ then
$ht(\lambda)=\alpha_1=\alpha_2$. We note that $detV(a)^k$ is a
semi-invariant of weight $(k,0)$. Hence
$(S_{\lambda}V_1)^{SL(V_1)}\otimes(S_{\lambda}V_2)^{SO(V_2)}$ is
generated by the semi-invariant $detV(a)^k$, so
$OSI(Q,\alpha)=\Bbbk[detV(a)]$. Finally if $\alpha_1<\alpha_2$ we
note that $det(-V(a)^tV(a))^k$ is a semi-invariant of weight
$(k,0)$ so
$(S_{\lambda}V_1)^{SL(V_1)}\otimes(S_{\lambda}V_2)^{SO(V_2)}$ is
generated by the semi-invariant $det(-V(a)^tV(a))^k$ and thus $OSI(Q,\alpha)=\Bbbk[det(-V(a)^tV(a))].$\\
 Now we prove the
induction step. Let $X=ORep(Q,\alpha)$
 and by theorem \ref{fc} we obtain
$$
OSI(Q,\alpha)=\big(\Bbbk[X]\big)^{SO(V)}=
$$

$$
\bigoplus_{{\lambda(a_1),\ldots,\lambda(a_{n})\in
\Lambda}}(S_{\lambda(a_1)}V_1)^{SL(V_1)}\otimes
(S_{\lambda(a_1)}V_2^*\otimes
S_{\lambda(a_2)}V_2)^{SL(V_2)}\otimes
$$
$$
\cdots\otimes (S_{\lambda(a_{n-1})}V_n^*\otimes
S_{\lambda(a_n)}V_n)^{SL(V_n)}\otimes
(S_{\lambda(a_n)}V_{n+1})^{SO(V_{n+1})} ,
$$
The proof of this theorem is the same of the proof of the theorem
\ref{tfse} up to when we have to consider $\alpha_{n+1}$. As in
the proof of theorem \ref{tfse} we can suppose
$\alpha_1\leq\cdots\leq \alpha_{n}$, otherwise, by induction, we
can reduce to a smaller quiver.\\
By lemma \ref{cl}, we have the generators
$det\,V(a_i)=det\,V(\sigma(a_i))$ if
$\alpha_i=\alpha_{i+1}$, $1\leq i\leq n-2$.\\
By proposition \ref{i2}, we obtain that there exist
$k_1,\ldots,k_{n}\in\mathbb{N}$ such that
$$
\lambda(a_{i})=(\overbrace{k_{i}+\cdots+k_1,\ldots,k_{i}+\cdots+k_1}^{\alpha_1},\ldots,\overbrace{k_{i},\ldots,k_{i}}^{\alpha_{i}-\alpha_{i-1}}),
$$
for every $i\in\{1,\ldots,n\}$.\\
Now, by proposition \ref{i3}, $\lambda(a_n)$ must have even rows.
If
 $\alpha_{n+1}=\alpha_j\leq \alpha_{j+1}\leq\cdots\leq \alpha_{n}$
 for some $j\in\{1,\ldots,n\}$ then $S_{\lambda(a_{n})}V_{n+1}^*=0$
 unless $ k_{n}+\cdots+k_{j+1}=0$, so
 $\lambda(a_{n})=\cdots=\lambda(a_{j+1})=\lambda(a_j)$. As before, by lemma \ref{cl}, we can consider the smaller
 quiver $Q^1:1\longrightarrow 2\cdots j\longrightarrow n+1\longrightarrow 2n-j+2\cdots 2n\longrightarrow 2n+1$  and then
 $$
OSI(Q,\alpha)\cong
(S_{\lambda(a_1)}V_1)^{SL(V_1)}\otimes\cdots\otimes(S_{\lambda(a_{j-1})}V_j^*\otimes
S_{\lambda(a_j)}V_j)^{SL(V_j)}\otimes
$$
\begin{equation}\label{02}
(S_{\lambda(a_{j-1})}V_n^*\otimes
S_{\lambda(a_j)}V_n)^{SL(V_n)}\otimes(S_{\lambda(a_j)}V_{n+1})^{SO(V_{n+1})},
\end{equation}
where
$$
\lambda(a_{j})=(\overbrace{k_{j}+\cdots+k_1,\ldots,k_{j}+\cdots+k_1}^{\alpha_1},\ldots,\overbrace{k_{j},\ldots,k_{j}}^{\alpha_{n+1}-\alpha_{j-1}}),
$$
and $k_{j}+\cdots+k_1,\ldots,k_j$ have to be even otherwise, by
proposition \ref{i3}, $(S_{\lambda(a_j)}V_{n+1})^{SO(V_{n+1})}=0$.
Hence $k_l$ has to be even for every $l\in\{1,\ldots,j\}$.
 Now to complete the proof it's enough to find the generators of the algebra (\ref{02}) for $\alpha_{n+1}=\alpha_j\leq \alpha_{j+1}\leq\cdots\leq
 \alpha_{n}$.\\
By proposition \ref{i2}, for every $l\in \{1,\ldots,j\}$,
$(S_{\lambda(a_{l-1})}V_l^*\otimes S_{\lambda(a_l)}V_l)^{SL(V_l)}$
is generated by a semi-invariant of weight
$(0,\ldots,0,k_l,0,\dots,0)$ where $k_l\in 2 \mathbb{N}$, is
$l$-th component. We note that
$$V(a_{2n-l+1})\cdots V(a_l)\in OSI(Q,\alpha)\mapsto
(det(V(a_{2n-l+1})\cdots V(a_l)))^{\frac{k_l}{2}}$$ is a
semi-invariant of weight
 $(0,\ldots,0,k_l,0,\dots,0)$, so it
generates $(S_{\lambda(a_{l-1})}V_l^*\otimes
S_{\lambda(a_l)}V_l)^{SL(V_l)}$. Since
$\lambda(a_l)=\lambda(a_{l-1})+(k_l)^{\alpha_l}$, then
$det(V(a_{2n-l+1})\cdots V(a_l))$ is a generator of
$OSI(Q,\alpha)$.\\
In the summand of $OSI(Q,\alpha)$ indexed by the families of
partitions in which
$\lambda(a_j)=(\overbrace{k_j,\ldots,k_j}^{\alpha_j=\alpha_{n+1}})$,
with $k_j\in 2\mathbb{N}$, we have that
$(S_{\lambda(a_{j})}V_j)^{SL(V_j)}\otimes
 (S_{\lambda(a_j)}V_{n+1})^{SO(V_{n+1})}$ is generated by
 a semi-invariant of weight
 $(0,\ldots,0,k_j,0,\dots,0,0)$ where $k_j$ is
  the $j$-th component and we note, as
 before, that $(det(V(a_{n})\cdots V(a_j)))^{k_j}$ is a
 semi-invariant of weight
 $(0,\ldots,0,k_j,0,\dots,0,0)$. Since
 $\lambda(a_j)=\lambda(a_{j-1})+(k_j)^{\alpha_j=\alpha_{n+1}}$,
 $det(V(a_{n})\cdots V(a_j))$ is a generator of $OSI(Q,\alpha)$.
\end{proof}
\begin{proof}[Proof of Theorem \ref{tfoij2}] Using the theorem
\ref{tfoe1}, the proof is similar to the proof of theorem
\ref{tfsij1}. \end{proof}


\begin{thebibliography}{99}
\bibitem{abw} K. Akin, D. Buchsbaum, J. Weyman, \textit{Schur
functors and Schur complexes}, Advances Math. 44 (1982), 207-277.
\bibitem{ass} I. Assem, D. Simson, A. Skowronski, \textit{Elements of the
Representation Theory of Associative Algebras}, volume 1, London
  Mathematical Society Students Texts 65, Cambridge University Press, 2006.
\bibitem{bgp} I. N. Bernstein, I. M. Gelfand, V. A. Ponomarev, \textit{Coxeter functors and Gabriel's theorem}, Uspekhi Mat. Nauk 28, no. 2(170) (1973), 19-33.
\bibitem{bmrrt} A. B. Buan, R. Marsh,
M. Reineke, I. Reiten, G. Todorov, \textit{Tilting theory and
cluster combinatorics}, Adv. Math. 204 (2006), 572-618.
\bibitem{dr} V. Dlab, C. M. Ringel, \textit{Indecomposable Representations of Graphs and Algebras}, Memoirs
Amer. Math. Soc. 173 (1976).

\bibitem{dw1} H. Derksen, J. Weyman, \textit{Semi-invariants of quivers and saturation for Littlewood-Richardson coefficients}, J. Amer. Math. Soc. 16 (2000),467-479.
\bibitem{dw2} H. Derksen, J. Weyman, \textit{Generalized
  quivers associated to reductive groups}, Colloq. Math. 94
  (2002), No. 2, 151-173.
  \bibitem{dz} M. Domokos, A. N. Zubkov, \textit{Semi-invariants of quivers as determinants}, Transform. Groups 6 (2001), no. 1, 9-24.  
  \bibitem{f} W. Fulton, \textit{Young tableaux, with
  applications to representation theory and geometry}, London
  Mathematical Society Student Texts 35, Cambridge University
  Press, 1997.
  \bibitem{fh} W. Fulton, J. Harris, \textit{Representation Theory; the first course},
Graduate Texts in Mathematics 129, Springer-Verlag, 1991.
\bibitem{fz1} S. Fomin, A. Zelevinsky, \textit{Cluster Algebras I:
Foundations}, J. Amer. Math. Soc. 15 (2002), 497-529.
\bibitem{fz2} S. Fomin, A. Zelevinsky, \textit{Cluster Algebras II: Finite type
classification}, Invent. Math. 154 (2003), 63-121.
\bibitem{iotw} K. Igusa, K. Orr, G. Todorov, and J. Weyman,
\textit{Cluster complexes via semi-invariants}, Compos. Math. 145 (2009), no. 4, 1001-1034.

\bibitem{l} A. A. Lopatin, \textit{Invariants of quivers under the action of classical groups}, J. Algebra 321 (2009), no. 4, 1079--1106.
 \bibitem{lz} A. A. Lopatin, A. N. Zubkov, \textit{Semi-invariants of mixed representations of quivers}, Transform. Groups 12 (2007), no. 2, 341-369.
 \bibitem{m} I. G. Macdonald, \textit{Symmetric functions and Hall polynomials},
second edition, with contributions by A. Zelevinsky, Oxford
Mathematical Monographs, Oxford University Press, New York, 1995.
\bibitem{mwz} P. Magyar, J. Weyman and A. Zelevinsky, \textit{Multiple
flag varieties of finite type}, Adv. Math. 141 (1999), no. 1,
97-118.

  \bibitem{p} C. Procesi, \textit{Lie Groups. An Approach through Invariants and Representations}, Universitext, Springer, New York, 2007.
\bibitem{r1} C. M, Ringel, \textit{Representations of K-species and
bimodules}, J. Algebra 41 (1976), 269-302.
\bibitem{s} A. Schofield, \textit{Semi-invariants of quivers}, J. London Math. Soc. (3) 65 (1992), 46-64.
\bibitem{sk} M. Sato, T. Kimura, \textit{A classification of irreducible prehomogeneous vector spaces and
their relative invariants}, Nagoya J. Math 65 (1977), 1-155.
\bibitem{sw} A. Skowronski, J. Weyman, \textit{The algebras
  of semi-invariants of quivers}, Trans. Groups 5 (2000), no. 4,
  361-402.
\end{thebibliography}
\end{document}